\documentclass[10pt,epsfig]{amsart}
\usepackage{amsbsy,amssymb,amscd,amsfonts,latexsym,amstext,delarray,
amsmath,amsthm, epsfig}

\input xypic

\newtheorem{The}{Theorem}[section]
\newtheorem{Cor}[The]{Corollary}
\newtheorem{Pro}[The]{Proposition}

\newtheorem{Lem}[The]{Lemma}

\numberwithin{equation}{section}

\def\proof{\vspace{2ex}\noindent{\bf Proof.} }
\def\endproof{\relax\ifmmode\expandafter\endproofmath\else
\unskip\nobreak\hfil\penalty50\hskip.75em\hbox{}\nobreak\hfil\bull
{\parfillskip= 0pt \finalhyphendemerits= 0 \bigbreak}\fi}
\def\endproofmath$${\eqno\bull$$\bigbreak}
\def\bull{\vbox{\hrule\hbox{\vrule\kern3pt\vbox{\kern6pt}\kern3pt\vrule}
\hrule}}

\newcommand{\no}{\noindent}

\newcommand{\ba}{\begin{eqnarray}}
\newcommand{\na}{\end{eqnarray}}
\newcommand{\scr}{\mathcal}

\newcommand{\bA}{\mathbb{A}}

\newcommand{\C}{\mathbb{C}}
\newcommand{\R}{\mathbb{R}}
\newcommand{\Z}{\mathbb{Z}}
\newcommand{\N}{\mathbb{N}}

\newcommand{\bP}{\mathbb{P}}

\renewcommand{\H}{\mathbb{H}}
\newcommand{\Q}{\mathbb{Q}}
\newcommand{\hL}{\mathbb{L}}
\newcommand{\fX}{{\mathfrak{X}}}

\newcommand{\D}{\Delta}
\newcommand{\cg}{{\mathfrak g}}
\newcommand{\gsl}{{\rm sl}}
\newcommand{\ie}{{\it i.e.\/}\ }
\newcommand{\eg}{{\it e.g.\/}\ }
\newcommand{\cf}{{\it cf.\/}\ }

\newcommand{\cA}{\scr{A}}

\newcommand{\mX}{\mathfrak{X}}

\newcommand{\T}{\scr{T}}

\newcommand{\sC}{\scr{C}}
\newcommand{\fC}{\mathfrak{C}}

\newcommand{\sP}{\scr{P}}

\newcommand{\sE}{\scr{E}}
\newcommand{\sK}{\scr{K}}

\newcommand{\sR}{\scr{R}}

\newcommand{\A}{\scr{A}}
\newcommand{\M}{\scr{M}}
\newcommand{\sO}{\scr{O}}
\newcommand{\sD}{\scr{D}}

\newcommand{\sH}{\scr{H}}
\newcommand{\Spec}{{\rm Spec}}
\newcommand{\Sp}{{\rm Spec}}
\newcommand{\Ker}{{\rm Ker}}

\newcommand{\Coker}{{\rm Coker}}
\newcommand{\ind}{{\rm ind}}
\newcommand{\SL}{{\rm SL}}

\newcommand{\Tr}{{\rm Tr}}
\newcommand{\Aut}{{\rm Aut}}

\title[Archimedean cohomology revisited]{Archimedean cohomology
revisited}
\author[Consani]{Caterina Consani$^{\text{\ddag}}$}
\author[Marcolli]{Matilde Marcolli$^{\text{\dag}}$}
\thanks{\noindent$^{\text{\ddag}}$Partially supported by NSERC
grants 72016789, 72024520}
\thanks{$^{\text{\dag}}$Partially supported by Humboldt Foundation
Sofja Kovalevskaja Award}
\address{C.~Consani: University of Toronto \ \ Canada}
\email{kc\@@math.toronto.edu}
\address{M.~Marcolli: Max--Planck Institut f\"ur Mathematik Bonn  \ \
Germany}
\email{marcolli\@@mpim-bonn.mpg.de}

\begin{document}

\maketitle

\section{Introduction}

\no C. Deninger produced a unified description of the local
factors at arithmetic infinity and at the finite places where the
local Frobenius acts semi-simply, in the form of a Ray--Singer
determinant of a ``logarithm of Frobenius'' $\Phi$ on an infinite
dimensional vector space (the {\em archimedean cohomology}
$H^\cdot_{ar}(X)$ at the archimedean places, \cf \cite{Den}). The
first author gave a cohomological interpretation of the space
$H^\cdot_{ar}(X)$, in terms of a double complex $K^{\cdot,\cdot}$
of real differential forms on a smooth projective algebraic
variety $X$ (over $\C$ or $\R$), with Tate-twists and suitable
cutoffs, together with an endomorphism $N$, which represents a
``logarithm of the local monodromy at arithmetic infinity''.
Moreover, in this theory the cohomology of the complex ${\rm
Cone}(N)^\cdot$ computes real Deligne cohomology of $X$ (\cf
\cite{KC}). The construction of \cite{KC} is motivated by a
dictionary of analogies between the geometry of the tubular
neighborhoods of the ``fibers at arithmetic infinity'' of an
arithmetic variety $X$ and the geometric theory of the limiting
mixed Hodge structure of a degeneration over a disk. Thus, the
formulation and notation used in \cite{KC} for the double complex
and archimedean cohomology mimics the definition, in the geometric
case, of a resolution of the complex of nearby cycles and its
cohomology(\cf \cite{Steen}).

\smallskip

\no In Section 2 and 3 we give an equivalent description of
Consani's double complex, which allows us
to investigate further the structure induced on the complex and
the archimedean cohomology by the operators $N$, $\Phi$, and the
Lefschetz operator $\hL$. In Section 4 we illustrate the analogies
between the complex and   
archimedean cohomology and a resolution of the complex of nearby cycles 
in the classical geometry of an analytic degeneration with normal 
crossings over a disk. 
In Section 5 we show that, using the Connes--Kreimer formalism 
of renormalization, we can identify the
endomorphism $N$ with the residue of a Fuchsian connection, in analogy to the
log of the monodromy in the geometric case. 
In Section 6 we recall Deninger's approach to the archimedean 
cohomology through an interpretation as global sections of a real
analytic Rees sheaf over $\R$.
In Section 7 we show how the action of the endomorphisms $N$ and 
$\hL$ and the Frobenius operator $\Phi$ define a noncommutative manifold 
(a spectral triple in the sense of Connes), where 
the algebra is related to the $\SL(2,\R)$ 
representation associated to the Lefschetz $\hL$, the Hilbert space is
obtained by considering Kernel and Cokernel of powers of $N$, and the log
of Frobenius $\Phi$ gives the Dirac operator. The archimedean
part of the Hasse-Weil L-function is obtained from a zeta function of
the spectral triple. In Section 8 we outline some formal analogies between
the complex and cohomology at arithmetic infinity and the equivariant 
Floer cohomology of loop spaces considered in Givental's homological 
geometry of mirror symmetry.

\section{Cohomology at arithmetic infinity}

\smallskip

\no Let $X$ be a compact K\"ahler manifold of (complex) dimension
$n$. Consider the complex of $\C$-vector spaces
\begin{equation}\label{alghbar}
C^\cdot = \Omega_X^\cdot \otimes \C [U, U^{-1}]\otimes \C
[\hbar, \hbar^{-1}],
\end{equation}
where $\Omega_X^\cdot = \oplus_{p,q}\Omega_X^{p,q}$ is the complex
of global sections of the sheaves of $(p,q)$-forms on 
$X$, $\hbar$ and $U$ are formal
independent variables, with $U$ of degree two. Our choice of
notation wants to be suggestive of \cite{Giv}, in view of the
analogies illustrated in the last section of this paper. On
$C^\cdot$ we consider the total differential $\delta_C =d'_C + d_C
''$, where $d'_C =\hbar \, d$, with $d=\partial + \bar
\partial$ the usual de Rham differential
and $d_C ''= \sqrt{-1}(\bar\partial -\partial)$. The
hypercohomology $\H^\cdot(C^\cdot,\delta_C)$ is then simply given
by the infinite dimensional vector space $H^\cdot (X;\C)\otimes
\C[U,U^{-1}] \otimes \C[\hbar,\hbar^{-1}]$.

\smallskip

\no We also consider the
positive definite inner product
\begin{equation}\label{innprod}
\langle \alpha \otimes U^r\otimes \hbar^k, \eta \otimes U^s \otimes \hbar^t
\rangle := \langle \alpha, \eta \rangle \, \, \delta_{r,s}
\delta_{k,t},
\end{equation}
where $\langle \alpha, \eta \rangle$ denotes the Hodge inner
product on forms $\Omega_X^\cdot$, given by
\begin{equation}\label{innprodint}
 \langle \alpha, \eta \rangle := \int_X
\alpha \wedge * \,C(\bar\eta),
\end{equation}
with $C(\eta)=(\sqrt{-1}\, )^{p-q}$, for $\eta\in
\Omega^{p,q}_X$, and $\delta_{a,b}$ the Kronecker delta. 

\smallskip

\no We then introduce certain cutoffs on $C^\cdot$, which will
allow us to recover the complex at arithmetic infinity of \cite{KC}
from $C^\cdot$.

\smallskip

\no To fix notation, for fixed $p,q\in \Z_{\geq 0}$ with $m=p+q$, let
\begin{equation}\label{lambda-cut}
 \lambda(q,r):= \max\left\{0,2r+m,r+q\right\},
\end{equation}
where $2r+m$ is the total degree of the complex.
Let $\tilde\Lambda_{p,q}\subset \Z^2$ be the set of lattice points
satisfying
\begin{equation}\label{tildeLambda}
\tilde\Lambda_q=\{ (r,k)\in \Z^2: \,\, k\geq \lambda(q,r) \},
\end{equation}
for $\lambda(q,r)$ as in \eqref{lambda-cut}.

\smallskip

\no For fixed $(p,q)$ with $m=p+q$, let
$\fC^{m,2r}_{p,q}\subset C^\cdot$ be the complex linear
subspace given by the span of the elements of the form
\begin{equation}\label{sCm2rpq-elts}
\alpha \otimes U^r \otimes \hbar^k,
\end{equation}
where $\alpha \in \Omega^{p,q}_X$ and $(r,k)\in
\tilde\Lambda_q$. We regard $\fC^{m,*}_{p,q}$ as a $2\Z$-graded
complex vector space.

\smallskip

\no Let $\fC^\cdot$ be the direct sum of all the
$\fC^{m,*}_{p,q}$, for varying $(p,q)$. We regard it as a $\Z$-graded
complex vector space with total degree $2r+m$.

\smallskip

\no In the cutoff \eqref{lambda-cut}, while the integer $2r+m$ is
just the total degree in $\Omega_X^\cdot \otimes \C[U,U^{-1}]$,
the constraint $k \geq r+q$ can be explained in terms of the Hodge
filtration.

\smallskip

\no Let $\gamma^\cdot =F^\cdot \cap \bar  F^\cdot$, where $F^\cdot$
and $\bar F^\cdot$ are the Hodge filtrations
\begin{equation}\label{Fprime}
F^p \Omega^m_X:=
\displaystyle{\bigoplus_{\substack{p'+q=m \\p'\geq p}}}
\,\, \Omega^{p',q}_X,
\end{equation}
\begin{equation}\label{Fprimeprime}
\bar F^q \Omega^m_X :=
\displaystyle{\bigoplus_{\substack{p+q''=m \\q''\geq q}}} \,\,
\Omega^{p,q''}_X.
\end{equation}

\smallskip

\no The condition defining $\fC^\cdot$ can be rephrased
in the following way.

\begin{Lem}\label{sC-Ffiltr}
The complex $\fC^\cdot$  has an equivalent
description as $\fC^i=\bigoplus_{i=m+2r} \fC^{m,2r}$, with
\begin{equation}\label{cutoffgammaC}
\fC^{m,2r} = \displaystyle{\bigoplus_{\substack{p+q=m \\k\geq
\max\{ 0, 2r+m\}}}} \,\,\, \left(F^{m+r-k} \, \Omega^m_X
\right) \otimes U^r \otimes \hbar^k .
\end{equation}
with the filtration $F^\cdot$ as in \eqref{Fprime}.
\end{Lem}

\proof This follows immediately by
\begin{equation}\label{primeFmrk}
 F^{m+r-k} \Omega^m_X =
\displaystyle{\bigoplus_{\substack{p+q=m \\ k\geq r+q }}} \,\,
\Omega^{p,q}_X.
\end{equation}
\endproof

\smallskip

\no Let $c$ denote the complex conjugation operator acting on
complex differential forms. We set $\T^\cdot :=
(\fC^\cdot)^{c=id}$. This is the real complex 
$$ (C^\cdot)^{c=id}=\Omega_{X,\R}^\cdot \otimes \R[U,U^{-1}]
\otimes \R[\hbar,\hbar^{-1}]. $$ 
Here $\Omega^m_{X,\R}$ is the $\R$-vector space of 
{\em real differential forms} of degree $m$,  
spanned by forms $\alpha=\xi+\bar\xi$, with $\xi \in
\Omega^{p,q}_X$ and such that $p+q=m$, namely
\begin{equation}\label{Omegam}
\Omega^m_{X,\R} =\bigoplus_{p+q = m} (\Omega^{p,q}_X+
\Omega^{q,p}_X).
\end{equation}
We have then the following equivalent description of
$\T^\cdot$.

\begin{Lem}\label{gammafiltr}
The complex $\T^\cdot=(\fC^\cdot)^{c=id}$ has the equivalent
description $\T^i = \bigoplus_{i=m+2r} \T^{m,2r}$ with
\begin{equation}\label{cutoffgamma}
\T^{m,2r} = \displaystyle{\bigoplus_{\substack{p+q=m \\k\geq
\max\{ 0, 2r+m\}}}} \,\,\, \left(\gamma^{m+r-k} \, \Omega^m_X
\right) \otimes U^r \otimes \hbar^k ,
\end{equation}
with the filtration $\gamma^\cdot =F^\cdot \cap \bar F^\cdot$.
\end{Lem}

\proof  We have
$$ F^{m+r-k} \Omega^m_X =
\displaystyle{\bigoplus_{\substack{p+q=m \\ k\geq r+q }}} \,\,
\Omega^{p,q}_X $$
$$ \bar F^{m+r-k} \Omega^m_X =
\displaystyle{\bigoplus_{\substack{p+q=m \\ k\geq r+p }}} \,\,
\Omega^{p,q}_X, $$
hence one obtains $$ \gamma^{m+r-k} \, \Omega^m_X
= \displaystyle{\bigoplus_{\substack{p+q=m \\ k\geq
\frac{|p-q|+2r+m}{2} }}} \,\, \Omega^{p,q}_X, $$
where $(|p-q|+2r+m)/2= r+ \max\{p,q\}$.

\endproof

\smallskip

\no Notice that the inner product \eqref{innprodint} is real
valued on real forms, hence it induces an
inner product on $\T^\cdot$.

\smallskip

\no Lemma \ref{gammafiltr} suggests the following convenient
description of $\T^\cdot$, which we shall use in the following.

\smallskip

\no For fixed $(p,q)$ with $m=p+q$, let
\begin{equation}\label{kappa}
 \kappa(p,q,r):= \max\left\{0,2r+m,\frac{|p-q|+2r+m}{2}\right\}.
\end{equation}
Then let
$\Lambda_{p,q}\subset \Z^2$ be the set
\begin{equation}\label{Lambda-pq}
 \Lambda_{p,q}=\{ (r,k)\in \Z^2 : \,\, k\geq \kappa(p,q,r) \},
\end{equation}
with $\kappa(p,q,r)$ as in \eqref{kappa}.

\smallskip

\begin{Pro}\label{TcomplCcompl}
The elements of the $\Z$-graded real vector space $\T^\cdot$ are
linear combinations of elements of the form
\begin{equation}\label{Tm2rpq-elts}
\alpha \otimes U^r \otimes \hbar^k,
\end{equation}
with $\alpha \in (\Omega^{p,q}_X+ \Omega^{q,p}_X)$,
$\alpha=\xi+\bar\xi$, for some $(p,q)$ with $p+q=m$, and $(r,k)$
in the corresponding $\Lambda_{p,q}$.
\end{Pro}

\proof By Lemma \ref{sC-Ffiltr}, we have seen that the cutoff
$\lambda(q,r)$ of \eqref{lambda-cut} corresponds to the Hodge
filtration $F^\cdot$, while Lemma \ref{gammafiltr} shows that, 
when we impose $c=id$ we can describe $(\fC^\cdot)^{c=id}$ in terms of the
$\gamma^\cdot$-filtration as in \eqref{cutoffgamma}. For fixed $(p,q)$
with $p+q=m$, this corresponds to the fact that
\begin{equation}\label{DefTmrpq}
 \T^{m,2r}_{p,q} := (\fC^{m,2r}_{p,q} \oplus \fC^{m,2r}_{q,p})^{c=id}
\end{equation}
is the real vector space generated by elements of the form
\eqref{Tm2rpq-elts}, where $\alpha=\xi+\bar\xi\in (\Omega^{p,q}_X+
\Omega^{q,p}_X)$ is a real form and the indices $(r,k)$ satisfy
the conditions $k\geq 0$, $k\geq 2r+m$ and $k\geq r+ \max\{p,q\}$.
Equivalently, $(r,k)\in \Lambda_{p,q}$. Notice that, since
$\kappa(p,q,r)=\kappa(q,p,r)$, we have
$\Lambda_{p,q}=\Lambda_{q,p}$. In fact $\Lambda_{p,q}$ depends on
$(p,q)$ only through $|p-q|$ and $m=p+q$.

\endproof

\smallskip

\begin{figure}
\begin{center}
\epsfig{file=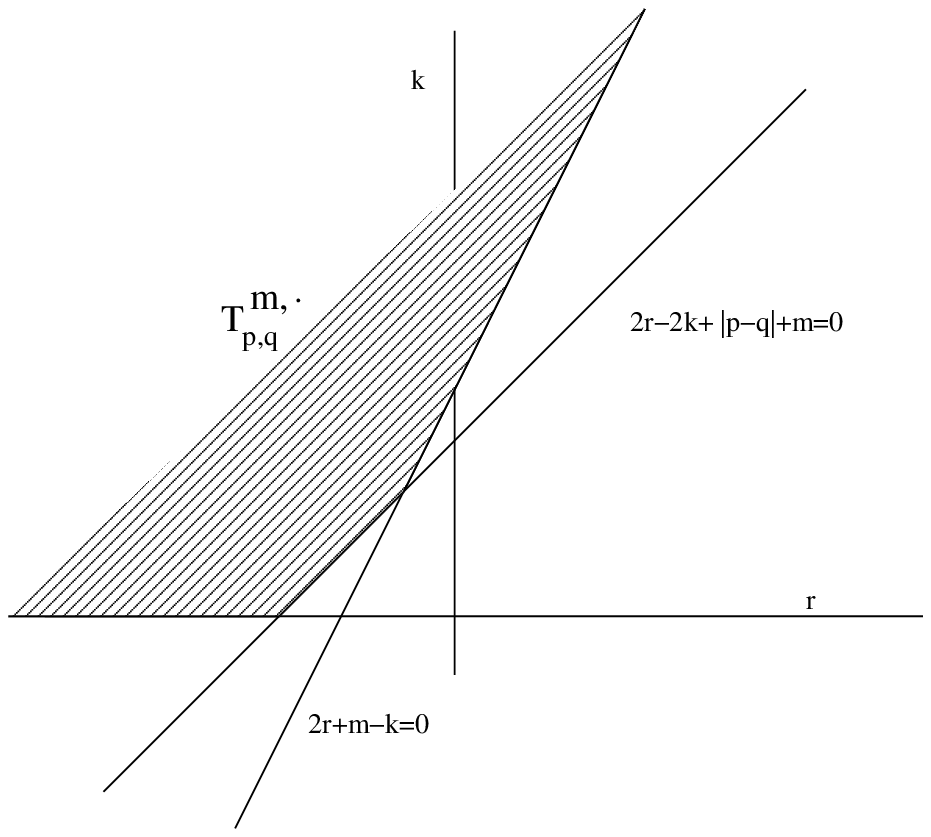} \caption{The region
$\Lambda_{p,q}\subset\Z^2$ defining $\T^{m,\cdot}_{p,q}$
\label{Figkappa}}
\end{center}
\end{figure}

\no Figure \ref{Figkappa} describes, for fixed values of $m$ and
$|p-q|$, the effect of the cutoff \eqref{kappa} on the varying indices
$(r,k)$. Namely, for fixed $(p,q)$ with $p+q=m$, the region
$\Lambda_{p,q}\subset \Z^2$ defined in \eqref{Lambda-pq} is the
shaded region in Figure \ref{Figkappa}.
The graph of the function $k=\kappa(p,q,r)$ of \eqref{kappa} is the
boundary of the shaded region in the Figure.

\smallskip

\no By construction, the real vector space
$\T^{m,2r}_{p,q}$ is the linear span of the \eqref{Tm2rpq-elts} with
$(r,k)\in \Lambda_{p,q}$ and $\T^\cdot$ is the direct sum of all the
$\T^{m,*}_{p,q}$, for varying $(p,q)$, viewed as a $\Z$-graded real
vector space with total degree $2r+m$.
Namely, we can think of a single $\T^{m,*}_{p,q}$ as a
``slice'' of $\T^\cdot$ for fixed $(p,q)$, namely, for each $(p,q)$
there is a corresponding Figure \ref{Figkappa} and $\T^\cdot$ is
obtained when considering the union of all of them.

\smallskip

\no The differentials $d'_C$ and $d_C''$ induce corresponding
differentials $d'$ and $d''$ on $\T^\cdot$, where $d'=d'_C =\hbar d$
and $d''=P^\perp d_C''$, with $P^\perp$ the orthogonal
projection of $(C^\cdot)^{c=id}$ onto $\T^\cdot$. Notice that, since
$d'$ and $d''$ change the values of $(p,q)$, the differentials move
from one ``slice'' $\T^{m,*}_{p,q}\subset \T^\cdot$ to another.

\smallskip
\subsection{Operators}

\smallskip

\no In the formulation introduced above, we then obtain the very
simple description of the operators $N$ and $\Phi$ of \cite{KC} as
\begin{equation}\label{NPhi}
N=U \hbar  \ \ \ \ \  \Phi = -U \frac{\partial}{\partial U}.
\end{equation}
In \cite{KC} these represent, respectively, a logarithm of the local
monodromy and a logarithm of Frobenius at arithmetic infinity.
We consider the Hilbert space completion of $\T^\cdot$ in the inner product
induced by \eqref{innprod}. With a slight abuse of notation, we
still denote this Hilbert space by $\T^\cdot$. 
The linear operator $N$ satisfies
$N^*N=P_1$ and $NN^*=P_2$, where $P_1$ and $P_2$ denote, respectively,
the orthogonal projections onto the closed subspaces $\Ker(N)^\perp$
and $\Coker(N)^\perp$ of $\T^\cdot$.
The operator $N$ also has the property that $[N,d']=[N, d'']=0$. The
operator $\Phi$ is an unbounded, self adjoint operator with spectrum
$\Spec(\Phi)=\Z$. It also satisfies $[\Phi,d']=[\Phi, d'']=0$.

\smallskip

\no Notice that, unlike the differentials $d'$ and $d''$ that move
between different slices $\T^{m,*}_{p,q}$ of $\T^\cdot$, the monodromy
map $N$ does not change the values $(p,q)$.

\smallskip

\no This means that, in addition to the result of Corollary
4.4 of \cite{KC} on the ``{\em global}'' properties of
injectivity and surjectivity of the map $N:\T^\cdot \to \T^{\cdot+2}$,
we can also give an analogous ``{\em local}'' result describing the 
properties of the map
$N:\T^{m,*}_{p,q} \to \T^{m,*+2}_{p,q}$, restricted to an action 
on a fixed ``slice'' (\ie for fixed $p$ and $q$). In this case, 
we obtain the following result.

\begin{Pro}\label{Nkercoker}
The endomorphism $N:\T^{m,2r}_{p,q} \to \T^{m,2(r+1)}_{p,q}$ has
the following properties:
\begin{enumerate}
\item $N$ is surjective iff $r$ is in the range $r > -\max\{ p,q \}$
\item $N$ is injective iff $r$ is in the range $r < -\min\{ p,q \}$.
\end{enumerate}
\end{Pro}

\proof (1) For fixed $(p,q)$ with $p+q=m$, let $\Lambda_{p,q}\subset \Z^2$
denote the shaded region in Figure \ref{Figkappa}, as in \eqref{Lambda-pq}.
Let $Z_{p,q} \subset \Z^2$ denote the set of lattice
points $Z_{p,q}=\{ (r,k)\in \Z^2 : \,\,r> -\max\{p,q\} \}$.

\smallskip

\no The point $(-\max\{p,q\},0)\in  \Lambda_{p,q}$ is the intersection
point of the lines $k=0$ and $2r-2k+m+|p-q|=0$ in the boundary of
$\Lambda_{p,q}$. Thus, one sees that the only points in $(r,k)\in
\Lambda_{p,q}$ such that $(r-1,k-1)\notin \Lambda_{p,q}$ are those of
the form $(r,0)$ with $r\leq -\max\{ p,q \}$. This shows that
every point $(r,k)\in \Lambda_{p,q}\cap
Z_{p,q}$ has the property that $(r-1,k-1)\in \Lambda_{p,q}$, hence $N$
is surjective in the range $r> -\max\{ p,q \}$. It also shows
that, for every $r\leq -\max\{ p,q \}$, the point $(r,0)$
in $\Lambda_{p,q}$ is such that $(r-1,-1)\notin \Lambda_{p,q}$, so
that N cannot be surjective in the range $r\leq -\max\{ p,q \}$.

\smallskip

\no (2) The case of injectivity is proved similarly. Let
$$ W_{p,q}=\{ (r,k)\in \Z^2: \,\, r< -\min \{ p,q \} \}. $$
Notice that the only points $(r,k)\in \Lambda_{p,q}$ such that
$(r+1,k+1)\notin \Lambda_{p,q}$ are those on the boundary line
$k=2r+m$. The point $(-\min \{ p,q \},|p-q|)$ is the intersection
point of the lines $k=2r+m$ and $2r-2k+m+|p-q|=0$ in the boundary of
$\Lambda_{p,q}$. Thus, we see that every point $(r,k)\in
\Lambda_{p,q}\cap  W_{p,q}$ is such that
$(r+1,k+1)\in \Lambda_{p,q}$, and conversely, for all $r\geq -\min\{
p,q \}$ there exists a point $(r,k=2r+m)\in\Lambda_{p,q}$ such that
$(r+1,k+1)\notin \Lambda_{p,q}$, hence $N$ is injective in the range
$r< -\min\{ p,q \}$, while it cannot be injective for $r\geq -\min\{
p,q \}$.

\endproof

\medskip

\no The complex $(\T^\cdot,\delta)$ has another important structure,
given by the Lefschetz operator, which, together with the
polarization and the monodromy, endows $(\T^\cdot,\delta)$ with the
structure of a polarized
Hodge--Lefschetz module, in the sense of Deligne and Saito (\cf \cite{KC},
\cite{GNA}, \cite{Saito}). The Lefschetz endomorphism $\hL$ is
given by
\begin{equation}\label{LefschetzOp}
 \hL = (\cdot \wedge \omega)\,\, U^{-1}
\end{equation}
where $\omega$ is the canonical real closed (1,1)-form determined by the
K\"ahler structure.
The Lefschetz operator satisfies $[\hL,d']=[\hL,d'']=0$.

\smallskip

\no Notice that, unlike the monodromy operator $N$ that preserves the
``slices'' $\T^{m,*}_{p,q}$, the Lefschetz moves between different
slices, namely
$$ \hL : \T^{m,*}_{p,q} \to \T^{m+2,*-1}_{p+1,q+1}. $$

\smallskip
\subsection{Dualities}

\smallskip

\no There are two important duality maps on the complex $\T^\cdot$.
The first is defined on forms by
\begin{equation}\label{involution}
 S: \alpha \otimes U^r \otimes \hbar^{2r+m+\ell} \mapsto \alpha \otimes
U^{-(r+m)} \otimes \hbar^\ell,
\end{equation}
for $\alpha\in \Omega^m_X$,
and it induces, at the level of cohomology, the duality map of
Proposition 4.8 of \cite{KC}.
The map $S$ induces, in particular, the duality between kernel and
cokernel of the monodromy map in cohomology, described in \cite{KC}
and \cite{CM} in terms of powers of the monodromy (\cf
Proposition \ref{dual2} below).

\smallskip

\no The other duality is given by the map
\begin{equation}\label{tildeS}
\tilde S : \alpha \otimes U^r \otimes \hbar^k \mapsto C(*\alpha)
\otimes U^{r-(n-m)} \otimes \hbar^k.
\end{equation}

\begin{Pro}\label{dual2}
Let $S$ and $\tilde S$ be the maps defined in \eqref{involution} and
\eqref{tildeS}.
\begin{enumerate}
\item The map $S: \T^\cdot \to \T^\cdot$ is an involution, namely
$S^2=1$. It gives a collection of linear isomorphisms
$$ S= N^{-(2r+m)}: {\rm span}\{ \alpha \otimes U^r \otimes
\hbar^{2r+m+\ell} \} \to {\rm span} \{ \alpha \otimes
U^{-(r+m)} \otimes \hbar^\ell \}, $$
realized by powers $N^{-(2r+m)}$ of the monodromy.
\item The map $\tilde S : \T^\cdot \to \T^\cdot$ is an involution,
$\tilde S^2 =1$. The map induced by $\tilde S$ on the primitive part of the
cohomology, with respect to the Lefschetz decomposition, agrees (up to
a non-zero real constant) with the power $\hL^{n-m}$ of the Lefschetz
operator.
\end{enumerate}
\end{Pro}

\begin{proof}
(1) The result for $S$ follows directly from the definition, in fact,
we have
$$ S^2(\alpha \otimes U^r \otimes
\hbar^{2r+m+\ell}) = S (\alpha \otimes
U^{-(r+m)} \otimes \hbar^\ell) $$
$$ = \alpha \otimes U^{r+m-m} \otimes \hbar^{\ell+2(r+m)-m} =
\alpha \otimes U^r \otimes \hbar^{2r+m+\ell}. $$
This means that the duality $S$ preserve the ``slices''
$\T^{m,*}_{p,q}$ and on them it
can be identified with the symmetry of $\T^\cdot$ obtained by
reflection of the shaded area of Figure \ref{Figkappa} along the line
illustrated in Figure \ref{FigkappaS}. The elements of the form
$\alpha \otimes U^{-m/2} \otimes \hbar^{|p-q|}$, for $\alpha \in
\Omega_X^m$ are fixed by the involution $S$ since
$$ S(\alpha \otimes U^{-m/2} \otimes \hbar^{|p-q|})=
\alpha \otimes U^{\frac{m}{2} -m} \otimes
\hbar^{|p-q|+\frac{2m}{2}-m} =\alpha \otimes U^{-m/2} \otimes
\hbar^{|p-q|}. $$

\smallskip
\begin{figure}
\begin{center}
\epsfig{file=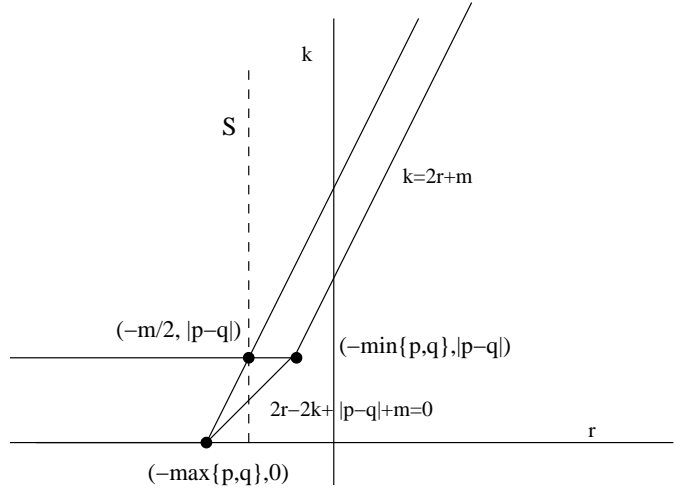} \caption{The duality $S$
\label{FigkappaS}}
\end{center}
\end{figure}

\smallskip

\no We prove (2). The map $\tilde S$ preserves $\T^\cdot$, since the
cutoffs described by the conditions $k\geq
0$, $2r+m -2k +|p-q|\leq 0$ and $k\geq 2r+m$ are preserved by mapping
$r\mapsto r-(n-m)$, $k\mapsto k$, $m\mapsto 2n-m$, $(p,q)\mapsto
(n-q,n-p)$. We have $\tilde S^2=1$, since
$$ \tilde S^2 (\alpha \otimes U^r \otimes \hbar^k)
= (\sqrt{-1})^{p-q}  \, \tilde S (*\alpha \otimes U^{r-(n-m)} \otimes
\hbar^k ) $$
$$ = (\sqrt{-1})^{p-q} (\sqrt{-1})^{n-q -(n-p)} \, *^2 \alpha \otimes
U^{r-(n-m)-(n-(2n-m))} \otimes \hbar^k $$
$$ = (\sqrt{-1})^{2(p-q)} (-1)^m\, \alpha \otimes U^{r}
\otimes \hbar^k, $$
where we used $*^2=(-1)^{m(2n-m)}=(-1)^m$.

\smallskip

\no Let $P^m(X)$ be the primitive part of the cohomology
$H^m(X,\C)$, with respect to the Lefschetz decomposition (\cf \cite{Wells}
\S V.6). Let $J$ be the operator induced on $P^{p,q}(X)$ by
$C(\eta)=(\sqrt{-1})^{p-q}\, \eta$. On the primitive cohomology one has
the identification (up to multiplication by a non-zero real constant)
\begin{equation}\label{starLef}
 \hL^{n-m} J \eta =  * \eta.
\end{equation}
In particular, \eqref{starLef} implies that the map $\tilde S$ agrees
(up to a normalization factor) with $\hL^{n-m}$, on the primitive
cohomology.

\end{proof}

\smallskip

\no Thus, we can think of the two dualities $S$ and $\tilde S$ as related,
respectively, to the action of $N^{2r+m}$ and $\hL^{n-m}$, \ie of
powers of the monodromy and Lefschetz.

\smallskip
\subsection{Representations}

\smallskip

\no The nilpotent endomorphisms $\hL$ and $N$ of $\T^\cdot$ introduced
above define two representations of $\SL(2,\R)$ on $\T^\cdot$ as
follows. With the notation
\begin{equation}\label{SL2Relts}
 \chi(\lambda) := \left(\begin{array}{cc} \lambda & 0 \\ 0 & \lambda^{-1}
\end{array} \right)\,\, \lambda \in \R^*, \ \
\ \ \  u(s):= \left(\begin{array}{cc} 1 & s \\ 0 & 1 \end{array}
\right)\,\, s\in \R,
\ \ \ \  w:= \left(\begin{array}{cc} 0 & 1 \\ -1 & 0 \end{array}
\right),
\end{equation}
we define $\sigma^L$ and $\sigma^R$ by
\begin{equation}\label{SL2L2}
\sigma^L(\chi(\lambda)) = \lambda^{-n+m}, \ \ \ \ \  \sigma^L(u(s))=
\exp(s\,\hL) ,
\ \ \ \ \ \sigma^L(w)= (\sqrt{-1})^n\,C\, \tilde S .
\end{equation}
\begin{equation}\label{SL2R2}
\sigma^R(\chi(\lambda)) = \lambda^{2r+m}, \ \ \ \ \  \sigma^R(u(s))=
\exp(s\,N) ,  \ \ \ \ \ \sigma^R(w)= C\, S .
\end{equation}
Here $C$ is the operator on forms $C(\eta)=(\sqrt{-1}\, )^{p-q}$
for $\eta\in \Omega^{p,q}_X$, and $S$ and $\tilde S$ the
dualities on $\T^\cdot$, as in Proposition \ref{dual2}. The
results of \cite{GNA}, \cite{KC} and \cite{CM} yield the
following.

\begin{Pro}\label{reps}
The operators \eqref{SL2L2} and \eqref{SL2R2} on $\T^\cdot$ define a
representation $\sigma=(\sigma^L,\sigma^R): \SL(2,\R)\times \SL(2,\R)\to {\rm
Aut}(\T^\cdot)$. In the representation $\sigma^L$,
the group $\SL(2,\R)$ acts by bounded operators on the completion of
$\T^\cdot$ in the inner product \eqref{innprod}.
Both $\SL(2,\R)$ actions commute with the Laplacian
$\Box= \delta \delta^* + \delta^* \delta$, hence they define induced
representations on the cohomology $\H^\cdot(\T^\cdot,\delta)$.
\end{Pro}

\proof For completeness, we give here a simple proof of the
proposition. In order to show that we have representations of
$\SL(2,\R)$ it is sufficient (\cite{LangSL2} \S XI.2) to check that
\eqref{SL2L2} and \eqref{SL2R2} satisfy the relations
\begin{equation}\label{LangRel}
\begin{array}{rl}
\sigma(w)^2 = & \sigma(\chi(-1)) \\[2mm]
\sigma(\chi(\lambda))\sigma(u(s))\sigma(\chi(\lambda^{-1})) = &
\sigma(u(s\lambda^2)) \end{array}
\end{equation}
We show it first for $\sigma^R$. We show that we have
$\sigma^R(w)=(-1)^m$, as in \cite{CM}. This follows directly from
the fact that $S^2=1$, since $\sigma^R(w)^2 = C S C S =
(\sqrt{-1})^{2(p-q)}= (-1)^m$. Thus, we have $\sigma^R(w)^2
=(-1)^m =\lambda^{2r+m}|_{\lambda=-1}$ and the first relation of
\eqref{LangRel} is satisfied. To check the second relation notice
that, on an element $\alpha\otimes U^r \otimes \hbar^k$ with
$\alpha\in \Omega^m_X$ we have
$$ \sigma^R(\chi(\lambda))\,\sigma^R(u(s))\,\sigma^R(\chi(\lambda^{-1}))
\,\, \alpha\otimes U^r \otimes \hbar^k = $$
$$ \sigma^R(\chi(\lambda)) \left( 1 + s N + \frac{s^2}{2} N^2 + \cdots
\right) \, \lambda^{-(2r+m)} \,\, \alpha\otimes U^r \otimes \hbar^k = $$
$$ \left( 1 + \lambda^{2(r+1)+m} s N \lambda^{-(2r+m)} +
\lambda^{2(r+2)+m} \frac{s^2}{2} N^2 \lambda^{-(2r+m)} + \cdots \right
)\, \alpha\otimes U^r \otimes \hbar^k =
 \exp( s\lambda^2\, N) \,\, \alpha\otimes U^r \otimes \hbar^k, $$
hence the second relation is satisfied.

\smallskip

\no We show that $\sigma^L$ also
satisfies the relations \eqref{LangRel}. Again, we first show that
$\sigma^L(w)=(-1)^{n+m}$, as in \cite{KC}. We have
$$ C\tilde S (\alpha\otimes U^r \otimes \hbar^k)=
(\sqrt{-1})^{n-q-n+p} (\sqrt{-1})^{p-q} *\alpha \otimes U^{r-(n-m)}
\otimes \hbar^k, $$
hence
$$ \sigma^L(w)^2\, (\alpha\otimes U^r \otimes \hbar^k) =
(-1)^n (-1)^m\, C\tilde S (*\alpha \otimes U^{r-(n-m)}
\otimes \hbar^k) $$
$$ = (-1)^n (-1)^m (-1)^m (\sqrt{-1})^{p-q} (\sqrt{-1})^{n-q-n+p}
\alpha \otimes U^r \otimes \hbar^k = (-1)^n (-1)^m \alpha \otimes U^r
\otimes \hbar^k $$
where, in the left-hand side, we used $*^2 =(-1)^m$.

\smallskip

\no Thus, we have $\sigma^L(w)^2= (-1)^{n+m}=
\lambda^{m-n}|_{\lambda=-1}$. Moreover, we have
$$ \sigma^L(\chi(\lambda))\,\sigma^L(u(s))\,\sigma^L(\chi(\lambda^{-1}))
\,\, \alpha\otimes U^r \otimes \hbar^k = $$
$$ \sigma^L(\chi(\lambda)) \left( 1 + s \hL + \frac{s^2}{2} \hL^2 + \cdots
\right) \, \lambda^{n-m} \,\, \alpha\otimes U^r \otimes \hbar^k = $$
$$ \left( 1 + \lambda^{-n+m+2} s \hL \lambda^{n-m} +
\lambda^{-n+m+4} \frac{s^2}{2} \hL^2 \lambda^{n-m} + \cdots \right
)\, \alpha\otimes U^r \otimes \hbar^k =
 \exp( s\lambda^2\, \hL) \,\, \alpha\otimes U^r \otimes \hbar^k, $$
hence the second relation is also satisfied.

\smallskip

\no The fact that $\sigma^L(\chi(\lambda))$ is a bounded operator in the
inner product induced by \eqref{innprod}, while for $\lambda\neq \pm
1$ the operators $\sigma^R(\chi(\lambda))$ are unbounded is clear from
the fact that the index $2r+m$ ranges over all of $\Z$, while $-n \leq
m-n \leq n$. For the fact that $[\Box,\sigma^L]=[\Box,\sigma^R]=0$ we
refer to \cite{KC}.

\endproof

\smallskip
\subsection{Ring of differential operators}

\smallskip

\no Let $\sD$ denote the algebra of differential operators on a
1-dimensional complex torus $T_\C$, generated by the operators
$Q=e^z$ and $P=\frac{\partial}{\partial z}$ satisfying the
commutation relation
\begin{equation}\label{PQ}
 PQ - QP =  Q.
\end{equation}
Let $\sR$ be the ring of functions defining the coefficients of
the differential operators in $\sD$. This is a subring of the ring
of functions on $\C^*$. For $\sR=\C[Q]$ we obtain $\sD=\C[P,Q]/(PQ
- QP =  Q)$.

\smallskip

\no Since the operators $N$ and $\Phi$ satisfy the commutation
relation $[\Phi,N]=-N$ and the operators $\hL$ and $\Phi$ satisfy
the commutation relation $[\Phi, \hL]= \hL$, the pairs of
operators $(N,-\Phi)$ and $(\hL,\Phi)$ define actions $\pi^R$ and
$\pi^L$ of $\sD$ on the complex $\fC^\cdot$ and on its cohomology,
by setting
\begin{equation}\label{actionD}
\begin{array}{ll}
\pi^L(P)=\Phi & \pi^L(Q)=\hL \\[3mm]
\pi^R(P)=-\Phi & \pi^R(Q)=N.
\end{array}
\end{equation}
There is an induced action of the ring $\sD=\R[P,Q]/(PQ
- QP =  Q)$ on the complex $\T^\cdot$ and on its cohomology.

\smallskip
\subsection{Weil--Deligne group at arithmetic infinity}

\smallskip

\no On $\T^\cdot$ consider the ``Frobenius flow''
\begin{equation}\label{Fflow}
F_t =e^{t\Phi}, \ \ \ \forall t\in \R,
\end{equation}
generated by the operator $\Phi$. We write
$F=F_1$. This satisfies
\begin{equation}\label{FN}
F \, N \, F^{-1}= e^{-1}\, N.
\end{equation}
Thus, the operators $F$ and $N$ can be thought of as defining an analog at
arithmetic infinity of the Weil--Deligne group ${\mathbb G}_a \rtimes
W_K$, which acts on the finite dimensional vector space associated to
the \'etale cohomology of the geometric generic fiber of a local
geometric degeneration for $K$ a non-archimedean
local field. In fact, in that case, the action of the 
Frobenius $\varphi \in W_K$ on $\mathbb G_a$ is given by
\begin{equation}\label{actionNF}
 \varphi \, x\, \varphi^{-1} = q^{-1} \, x,
\end{equation}
where $q$ is the cardinality of the residue field. The formal
replacement of $q$ by $e$ and of $\varphi$ by
$F$ determines \eqref{FN} from \eqref{actionNF}. In the
archimedean case, this ``Weil--Deligne group'' acts directly at the
level of the complex, not just on the cohomology.

\section{Archimedean cohomology}

\smallskip

\no We now describe the relation between the complex
$(\T^\cdot,\delta=d'+d'')$ defined in the first section
 and the cohomology theory
at arithmetic infinity developed in \cite{KC}.

\smallskip

\no On a smooth projective algebraic variety $X$ of dimension $n$
over $\C$ or $\R$, the complex of Tate-twisted real differential
forms introduced in (4.1) of \cite{KC} is defined as
\begin{equation}\label{Kcomplex}
K^{i,j,k} = \begin{cases}
\displaystyle{\bigoplus_{\substack{p+q=j+n\\|p-q|\le
2k-i}}}(\Omega^{p,q}_X\oplus\Omega^{q,p}_X)_\R\, \otimes_\R
\R\left(\frac{n+j-i}{2}\right)
&\text{if $j+n-i\equiv~0(2)$, $k\ge\text{max}(0,i)$} \\
0 &\text{otherwise,}
\end{cases}
\end{equation}
for $i,j,k\in \Z$. Here
$\R(r)$ denotes the real Hodge structure $\R(r) :=
(2\pi\sqrt{-1})^r \R$, and the differentials $d'$ and $d''$ are
given by
\begin{equation}\label{d1K}
d': K^{i,j,k} \to K^{i+1,j+1,k+1},\qquad d'(\alpha) = d(\alpha)
\end{equation}
\begin{equation}\label{d2K}
d'': K^{i,j,k} \to K^{i+1,j+1,k},\qquad d''(\alpha) =
\sqrt{-1}(\bar\partial -\partial)(\alpha)\quad(\text{projected onto
$K^{i+1,j+1,k}$}).
\end{equation}
The inner product on $\H^\cdot(K^\cdot,d'+d'')$ is
defined in terms of the bilinear form
\begin{equation}\label{hLinnerprod}
 Q(\alpha,\eta)=\int_X \hL^{n-m} \alpha \wedge J \bar\eta ,
\end{equation}
for $\eta,\alpha$ in the primitive part $P^m(X)$ of the de Rham
cohomology $H^m(X,\R)$ with respect to the Lefschetz
decomposition.

\smallskip

\no The relation between the complex $\T^\cdot$ and the total complex
$K^\cdot$ of \eqref{Kcomplex} is described by the following result,
which shows that the complex $\T^\cdot$ is identified with $K^\cdot$ 
after applying the simple change of variables
\begin{equation}\label{indicesrmk}
i=m+2r, \ \ \ \ r=-\frac{n+j-i}{2}, \ \ \ \ j=-n+m
\end{equation}
to the indices of $\T^\cdot$, taking fixed points
under complex conjugation ($c=id$) and replacing the variable $U$
by a Tate twist.

\begin{Pro}\label{TK}
Upon identifying the formal variable $U^{-1}$ with the Tate
twist given by multiplication by $2\pi \sqrt{-1}$, we obtain an
isomorphism of complexes
\begin{equation}\label{T=K}
(\T^\cdot |_{U=(2\pi \sqrt{-1})^{-1}},\delta) \cong
(K^\cdot,d'+d''),
\end{equation}
given by a reparameterization of the indices. Up to a
normalization factor, the inner product induced by \eqref{innprod}
on the cohomology $\H^\cdot(\T^\cdot |_{U=(2\pi
\sqrt{-1})^{-1}},\delta)$ agrees with the one defined by
\eqref{hLinnerprod} on $\H^\cdot(K^\cdot,d'+d'')$.
\end{Pro}

\begin{proof} We define a homomorphism
$$ I: \T^\cdot \to K^\cdot $$
as follows. For fixed $p,q$ with $p+q=m$, consider the region
$\Lambda_{p,q}\subset \Z^2$ as in \eqref{Lambda-pq}, with
$\kappa(p,q,r)$ as in \eqref{kappa}. For every $\alpha\in
(\Omega^{p,q}_X \oplus \Omega^{q,p}_X)_\R$ a real form,
$\alpha=\xi+\bar\xi$, with $\xi\in \Omega^{p,q}_X$, and for
every point $(r,k)\in \Lambda_{p,q}$, we have
$$ \alpha \otimes U^r \otimes \hbar^k \in \T^{m,2r}_{p,q}, $$
and, by Proposition \ref{TcomplCcompl},
every element of $\T^\cdot$ is a linear combination of elements of
this form, for varying $(p,q)$ and corresponding $(r,k)\in
\Lambda_{p,q}$.

\smallskip

\no We now define the map $I$ in the following way. To an element
\begin{equation}\label{etaxi}
 (\xi+\bar\xi) \otimes U^r \otimes \hbar^k,
\end{equation}
with $\xi\in \Omega^{p,q}_X$, with $p+q=m$, and with $k\geq
\kappa(p,q,r)$, we assign an element $I(\eta)\in K^{i,j,k}$, with
the same index $k$ and with
\begin{equation}\label{indicesrmk1}
i=m+2r, \ \ \text{ and }  \ \ \ j=-n+m,
\end{equation}
by setting
\begin{equation}\label{Ieta}
I(\eta) = (2\pi \sqrt{-1})^{-r}\,  (\xi+\bar\xi).
\end{equation}
In fact, for $(i,j)$ as in \eqref{indicesrmk1}, the index $r\in \Z$
can be written in the form
$$ r=-\frac{n+j-i}{2}, \ \ \  \text{ where } n+j-i = 0 \mod 2. $$
Thus, the element \eqref{Ieta} can be written as
\begin{equation}\label{Ieta2}
 (2\pi \sqrt{-1})^{\frac{n+j-i}{2}} \, (\xi+\bar\xi).
\end{equation}
By the definition \eqref{Kcomplex}, to check that this is an element
in $K^{i,j,k}$, it is sufficient to verify that $p+q=j+n$,
and that the conditions $|p-q|\le 2k-i$ and $k\ge \max\{ 0, i\}$ are
satisfied. Since $j=-n+m$ and $p+q=m$ we have $p+q=j+n$. The index $k$
is the same as in \eqref{etaxi}, hence it satisfies $k\geq
\kappa(p,q,r)$. This means that $k\geq 0$ and that $k\geq 2r+m=i$, so
that $k\ge \max\{ 0, i\}$ satisfied. Since $k\geq
\kappa(p,q,r)$ also implies $k\geq (|p-q|+2r+m)/2$, which, by $i=2r+m$
is the condition $|p-q|\le 2k-i$. It is clear that the map $I$ defined
this way is injective.

\smallskip

\no Thus, we have shown that, for every real form $\alpha \in
(\Omega^{p,q}_X \oplus \Omega^{q,p}_X)_\R$ and for every
lattice point $(r,k) \in \Lambda_{p,q}$ we have a unique
corresponding element in the complex $K^{i,j,k}$.

\smallskip

\no The map $I$ is also surjective, hence a linear isomorphism. In
fact, every element in $K^\cdot$ is a linear combination of
elements of the form \eqref{Ieta2} in $K^{i,j,k}$, for $\xi\in
\Omega^{p,q}_X$, and indices $(i,j,k)\in \Z^3$ satisfying
$n+j-i = 0$ mod $2$, $p+q=j+n$, $k\geq \max\{ i,0 \}$ and
$|p-q|\le 2k-i$. It is sufficient to show that, for any such
element, there exists a point $(r,k)\in \Lambda_{p,q}$ such that
$$ I((\xi+\bar\xi)\otimes U^r\otimes \hbar^k)= (2\pi
\sqrt{-1})^{\frac{n+j-i}{2}} \, (\xi+\bar\xi) \in K^{i,j,k}. $$
This is achieved by taking the point
\begin{equation}\label{point}
 \left( r=-\frac{n+j-i}{2}, k \right).
\end{equation}
Since $n+j-i = 0$ mod $2$ this point is in $\Z^2$, and since under the
change of variables \eqref{indicesrmk1} the conditions $k\geq \max\{
i,0 \}$ and $|p-q|\le 2k-i$ are equivalent to the condition $k\geq
\kappa(p,q,r)$ of \eqref{kappa}, the point \eqref{point} is in
$\Lambda_{p,q}$.

\smallskip

\no The map $I$ is compatible with the differentials, namely
$$ I(\delta\eta)=(d'+d'') I(\eta), $$
where on the left hand side $\delta=d'_C + P^\perp d_C''$ is the
differential on $\T^\cdot$ and on the right hand side $d'+d''$ is as
in \eqref{d1K} \eqref{d2K}. To see this, first notice that the
differential $d'$ of \eqref{d1K} satisfies
$$ d' = I d_C' I^{-1} =I \hbar d I^{-1}: K^{2r+m,-n+m, k} \to
K^{2r+(m+1), -n + (m+1), k+1}. $$
The analogous statement $d'' = I P^\perp d_C'' I^{-1}$ for the differential
$d''$ of \eqref{d2K} also involves the fact that the orthogonal
projection onto $K^{i+1,j+1,k}$ in \eqref{d2K}, induced by the inner
product \eqref{hLinnerprod} agrees with the corresponding orthogonal
projection $P^\perp$ in $\T^\cdot$ induced by the inner product
\eqref{innprod}.

\smallskip

\no The identification \eqref{starLef} implies that the inner
product \eqref{hLinnerprod} on $\H^\cdot(K^\cdot,d'+d'')$
considered in \cite{KC} and the inner product induced by
\eqref{innprod} agree up to a normalization factor.

\end{proof}

\smallskip

\no In particular, the `weight type' condition
$|a-b|\le 2k-i$ on the real forms in \eqref{Kcomplex} describes, as in
\eqref{Fprime}
\eqref{Fprimeprime}, the
filtration $\gamma^\cdot := F^\cdot\cap \bar F^\cdot$ on
the complex of real differential forms  on $X$.
It follows that the complex $K^{i,j,k}$ has a real analytic type,
even when $X$ defined over $\C$ does not have a real structure.

\smallskip

\no To the abelian group $K^{i,j,k}$ we assign the weight: $
-n-j+i\in\Z$. Keeping in mind that $\R(\frac{n+j-i}{2})$ is the
real Hodge structure of rank one and pure bi-degree
$(-\frac{n+j-i}{2},-\frac{n+j-i}{2})$, we obtain the following
description in terms of the filtration $\gamma^\cdot$:
\begin{equation}\label{gammaKijk}
K^{i,j,k} = \gamma^{\frac{n+j+i}{2}-k}\Omega^{n+j}_X
\otimes_{\R} \R(\frac{n+j-i}{2}) =
\gamma^{\frac{n+j+i}{2}-k}\Omega^{n+j}_X \otimes_{\R}
\gamma^{-\frac{n+j-i}{2}}\R .
\end{equation}
When considering the tensor product of the two structures
one sees that the index of the $\gamma$-filtration on the product (\ie
on $K^{i,j,k}$) is $i-k$.

\smallskip
\subsection{Deligne cohomology}

\smallskip

\no By Proposition \ref{TK}, the complex $\T^\cdot$ is related to
the real Deligne cohomology $H_\sD^*(X,\R(r))$. These groups can
be computed via the Deligne complex $(C^*_\sD(r),d_\sD)$ (\cf
\cite{Burgos}, \cite{KC}). The relation of $(C^*_\sD(r),d_\sD)$ to
the complex $(\T^\cdot,\delta)$ is given by the following result
of \cite{KC} (Prop.~4.1), which, for convenience, we reformulate
here in our notation.

\begin{Pro}\label{coneDeligne}
For $N$ acting on $(\T^\cdot,\delta)$, consider the complex
$(Cone(N)^\cdot,D)$ with differential $D(\alpha,\beta) =
(\delta(\alpha), N(\beta)-\delta(\beta))$.
\begin{enumerate}
\item For $2r+m>0$, the map $N^{-(2r+m)}$ gives an isomorphism between
the cohomology group in $\H^\cdot(\Ker(N)^\cdot)$, which lies over the
point of coordinates $(2r,2r+m)$ in Figure \ref{Figkappa}, and
the cohomology group in $\H^\cdot(\Coker(N)^\cdot)$ that lies over the
point of coordinates $(-2(r+m),0)$ in Figure \ref{Figkappa}.
\item In the range $2r+m<-1$, the
cohomology $\H^\cdot(Cone(N)^\cdot,D)$ is identified with
$\H^{\cdot+1}(\Coker(N),\delta)$.
\item Upon identifying the variable $U^{-1}$ with the Tate twist by
$2\pi \sqrt{-1}$, and for a for fixed $r\in \Z_{\leq 0}$, we obtain quasi
isomorphic complexes
\begin{equation}\label{quasi-iso}
(Cone(N)^\cdot,D)|_{U^r=(2\pi \sqrt{-1})^{-r}} \simeq
(C^*_\sD(-r),d_\sD)
\end{equation}
\end{enumerate}
\end{Pro}

\smallskip
\subsection{Local factors}

\smallskip

\no The ``archimedean factor'' (\ie the local factor at arithmetic
infinity) $L_\kappa(H^m,s)$ is a product of powers of shifted
Gamma functions, with exponents and arguments that depend on the
Hodge structure on $H^m = H^m (X,\C)= \oplus_{p+q=m} H^{p,q}$.
More precisely, it is given by (\cf \cite{Serre})
\begin{equation}\label{factor}
 L_\kappa (H^m,s)= \begin{cases}
\prod_{p,q}\Gamma_\C(s-\text{min}(p,q))^{h^{p,q}}
&\kappa = \C  \\\\
\prod_{p<q}\Gamma_\C(s-p)^{h^{p,q}}\prod_p
\Gamma_\R(s-p)^{h^{p+}}\Gamma_\R(s-p+1)^{h^{p-}} & \kappa = \R,
\end{cases}
\end{equation}
where the $h^{p,q}$, with $p+q=m$, are the Hodge numbers,
$h^{p,\pm}$ is the dimension of the $\pm(-1)^p$-eigenspace of de
Rham conjugation on $H^{p,p}$, and
$$ \Gamma_\C(s) := (2\pi)^{-s}\Gamma(s) \ \ \ \ \
\Gamma_\R(s) :=2^{-1/2}\pi^{-s/2}\Gamma(s/2).
$$
It is shown in \cite{Den} that the local factor \eqref{factor} can
be computed as a Ray--Singer determinant
\begin{equation}\label{DenDet}
L_\kappa(H^m,s) =
\det_{\infty}\left(\frac{1}{2\pi}(s-\Phi)|_{H^m_{ar}(X)}
\right)^{-1},
\end{equation}
where the ``archimedean cohomology'' $H^m_{ar}(X)$ is an
infinite-dimensional real vector space, and the zeta regularized
determinant of
an unbounded self adjoint operator $T$ is defined as
$$\det_\infty(s-T)=\exp (-\frac{d}{dz} \zeta_T(s,z)|_{z=0}).$$ In
\cite{KC} (\cf \S 5) the archimedean cohomology is identified
with the ``inertia invariants''
\begin{equation}\label{archinert}
H^\cdot_{ar}(X)=\H^\cdot(K^\cdot,d'+d'')^{N=0},
\end{equation}
where $\H^\cdot(K^\cdot,d'+d'')$ is the hypercohomology of the
complex \eqref{Kcomplex} and $\H^\cdot(K^\cdot,d'+d'')^{N=0}$ is
the kernel of the map induced in hypercohomology by the monodromy
$N$. This follows the expectation that the fiber over arithmetic
infinity has semi-stable reduction.

\smallskip

\no At arithmetic infinity, the alternating product of the local
factors of $X$ can be described in terms of the operators
$\sigma^L(w)^2$ and $\Phi$ (\cf \cite{CM} par.~3.4). Let $\Phi_0$
denote the restriction of the operator $\Phi$ to
$\H^\cdot(K^\cdot)^{N=0}$. For $a$ a bounded operator on
$\H^\cdot(K^\cdot)^{N=0}$, let $\zeta_{a,\Phi_0}(s,z)$ denote the
two variable zeta function
\begin{equation}\label{zetasz}
\zeta_{a,\Phi_0}(s,z)=\sum_{\lambda\in \Spec(\Phi_0)} \Tr(a
\Pi_\lambda) (s-\lambda)^{-z},
\end{equation}
where $\Pi_\lambda$ are the spectral projections of $\Phi_0$. Let
$\det_{\infty,a,\Phi_0}$ denote the zeta regularized determinant
\begin{equation}\label{zetadeta}
\det_{\infty,a,\Phi_0}(s)= \exp\left(-\frac{d}{dz}
\zeta_{a,\Phi_0}(s,z) |_{z=0} \right).
\end{equation}

\begin{Pro}\label{zetaL}
The two variable zeta function $\zeta_{\sigma^L(w)^2,\Phi_0}(s,z)$
satisfies
\begin{equation}\label{zetadet}
\det_{\infty,\sigma^L(w)^2,\Phi_0/(2\pi)} \left( \frac{s}{2\pi}
\right)^{-1} = \prod_{m=0}^{2n} L_\C(H^m,s)^{(-1)^{m+n}}.
\end{equation}
\end{Pro}

\begin{proof}
The operator $\Phi_0$ has spectrum
$$ \Spec(\Phi_0)=\{ \lambda_{\ell,p,q}=\min\{ p,q \} - \ell : \ell\in
\Z_{\geq 0} \}, $$ with eigenspaces $E_{\ell,p,q}=H^{p,q}(X)
\otimes U^r \otimes \hbar^{2r+m}$, with $r=\ell-\min\{ p,q \}$ and
$m=p+q$. In fact, for $\H^\cdot(K)^{N=0}$ we have $2r+m=k$, and
$k\geq |p-q|$ so that $r \geq -\min\{ p,q \}$. The result then
follows as in \cite{CM} \S 3.3.

\end{proof}

\section{Archimedean cohomology and nearby cycles}

\no The definition of the complex \eqref{Kcomplex} was inspired by
an analogy with the resolution of the complex of nearby cycles
associated to an analytic degeneration with normal crossings over
a disk, \cf~\cite{Steen}. In this section we recall this classical
construction and we relate it to its archimedean counterpart by 
making the analogy more explicit.

\smallskip
\subsection{The complex of nearby cycles}

\smallskip

\no Let $\fX$ and $\D$ be complex analytic manifolds, of complex
dimensions $\dim \fX = n+1$ and $\dim \D = 1$, and let $f: \fX \to
\D$ be a flat, proper morphism with projective fibers. For $0\in
\D$, we write $Y = f^{-1}(0)$, $\fX^* = \fX \smallsetminus Y$, and
$\D^* = \D\smallsetminus \{ 0 \}$. We assume that the map $f$ is
smooth on $\fX^*$ and that $Y$ is a divisor with normal crossings
on $\fX$. Under these hypotheses, the relative de Rham complex of
sheaves of differential forms with logarithmic poles along $Y$ is
well defined and of the form
\[
\Omega_{\fX/\D}^m(\log Y):= \wedge^{m}\Omega^1_{\fX/\D}(\log Y),
\]
where we use the inclusion $f^*\Omega^1_{\D}(\log 0) \subset
\Omega^1_{\fX}(\log Y)$ to define
\[
\Omega^1_{\fX/\D}(\log Y) := \Omega^1_{\fX}(\log
Y)/f^*\Omega^1_{\D}(\log 0).
\]
Consider local coordinates $\{z_0,\ldots z_n\}$ at $P\in Y$, with
$t$ the local coordinate at $0\in \D$, so that $t\circ f =
z_0^{e_0}\cdots z_k^{e_k}$ for some $0\le k\le n$ and $e_i\in\N$.
The stalk $\Omega^1_{\fX/\D}(\log Y)_P$ is the ${\mathcal
O}_{\fX,P}$-module with generators $\{dz_0/z_0,\ldots,
dz_k/z_k,dz_{k+1},\ldots dz_n\}$ satisfying the relation
$\sum_{i=0}^k e_i dz_i/z_i =0$. Thus, $\Omega^1_{\fX/\D}(\log Y)$
is a locally free sheaf of rank $n = \dim \fX - 1$ endowed with
differential $d$ given by the composite
$$ d: {\mathcal O}_{\fX} \to \Omega^1_{\fX}
\hookrightarrow \Omega^1_{\fX}(\log Y) \to \Omega^1_{\fX/\D}(\log
Y),
$$ where $\Omega_{\fX}^1(\log Y)$ is the free sheaf of ${\mathcal
O}_{\fX}$-modules on the same generators of
$\Omega_{\fX/\D}^1(\log Y)$.

\smallskip

\no The relative hypercohomology sheaves $\mathbb R^m
f_*\Omega^\cdot_{\fX/\D}(\log Y)$ are locally free ${\mathcal
O}_{\D}$-modules of finite rank (\cf \cite{Steen}). This follows
from the fact that the restriction $f: \fX^* \to \D^*$, which is a
smooth fiber bundle, determines, for all $s\in \D^*$, a canonical
isomorphism of complexes
\[
\Omega^\cdot_{\fX/\D}(\log Y)\otimes_{{\mathcal O}_\fX}{\mathcal
O}_{\fX_s} \stackrel{\simeq}{\to}\Omega^\cdot_{\fX_s}
\]
on $\fX_s = f^{-1}(s)$. This implies that
$\Omega^\cdot_{\fX/\D}(\log Y)\otimes_{{\mathcal O}_\fX}{\mathcal
O}_{\fX_s}$ is a resolution of the constant sheaf $\C$ on
$\fX_{s}$, hence
\[
\mathbb H^m(\fX_s,\Omega^\cdot_{\fX/\D}(\log Y)\otimes_{{\mathcal
O}_{\fX}}{\mathcal O}_{\fX_s}) \simeq H^m(\fX_s,\C), \ \ \ \
\forall s\in \D^*,
\]
hence the complex dimension of ${\mathbb
H}^m(\fX_s,\Omega^\cdot_{\fX/\D}(\log Y)\otimes_{{\mathcal
O}_{\fX}}{\mathcal O}_{\fX_s})$ is a locally constant function. One
obtains
\begin{equation}\label{RmHm}
{\mathbb R}^m f_*\Omega^\cdot_{\fX/\D}(\log Y)\otimes_{{\mathcal
O}_{\D^*}}k(s) \simeq H^m(\fX_s,\C), \ \ \ \  \forall s\in \D^*,
\end{equation}
so that $\mathbb R^m f_*\Omega^\cdot_{\fX/\D}(\log Y)$ is a
locally free ${\mathcal O}_{\D^*}$-module of finite rank.
Moreover, if $\D$ is a small disk, there exists an isomorphism
\begin{equation}\label{isoX/S}
H^m (\tilde \fX^*,\C) \stackrel{\simeq}{\to}{\mathbb H}^m
(Y,\Omega^\cdot_{\fX/\D}(\log Y)\otimes_{{\mathcal
O}_{\fX}}{\mathcal O}_Y)
\end{equation}
where $\tilde \fX^* = \fX \times_{\D}\tilde \D^*$ and $\tilde \D^*
\to \D^*$ is the universal covering space of $\D^* =
\D\setminus\{0\}$, so that $\dim{\mathbb H}^m
(\fX_s,\Omega^\cdot_{\fX/\D}(\log Y)\otimes_{{\mathcal
O}_{\fX}}{\mathcal O}_{\fX_s})$ is locally constant on $\D$.

\smallskip

\no A stronger result (\cf \cite{Steen}, \S 9) shows that the
Gauss--Manin connection
\begin{equation}\label{GaussManin}
\nabla: {\mathbb R}^m f_*\Omega^\cdot_{\fX/\D}(\log Y) \to
\Omega^1_{\D}(\log 0) \otimes_{{\mathcal O}_{\D}} {\mathbb R}^m
f_*\Omega^\cdot_{\fX/\D}(\log Y)
\end{equation}
has logarithmic singularities at $0\in \D$ and in case of an
algebraic morphism $f$ admits an algebraic description. In
particular the residue of $\nabla$ at zero is a well defined
operator
\begin{equation}\label{Nresidue}
 N := \text{Res}_0(\nabla)
\end{equation}
in fact it is an endomorphism of $H^m (\tilde \fX^*,\C)$. The
eigenvalues of $\text{Res}_0(\nabla)$ are rational numbers
$\alpha$ with $0\le\alpha <1$. The monodromy transformation $T$
induces an automorphism $T_0$ of \eqref{RmHm} and it can be shown
that
\begin{equation}\label{expNres}
T_0 = \exp(-2\pi\sqrt{-1}\,\,\text{Res}_0(\nabla)).
\end{equation}
It follows that the eigenvalues of $T_0$ are roots of unity, so
that a power $T_0^d$ is unipotent. In fact, up to a base change
$\D \to \D$, $z \mapsto z^d$, one can assume that $T_0$ is
already unipotent, this means that the residue \eqref{Nresidue} is
nilpotent.

\smallskip

\no There are two important filtrations on the cohomology $H^m
(\tilde \fX^*,\C)$. One is the Hodge filtration $F^p H^m(\tilde
\fX^*,\C)$ determined by the isomorphism \eqref{isoX/S} and the
`naive filtration' on $\Omega^\cdot_{\fX/\D}(\log
Y)\otimes_{{\mathcal O}_{\fX}}{\mathcal O}_Y$. The other is the
Picard--Lefschetz filtration $L_\ell H^m (\tilde \fX^*,\C)$,
$\ell\in\Z$, which is a canonical, finite, increasing filtration
associated to the map $N$ and defined by induction. The properties
and behavior of this filtration are a priori rather mysterious as
its definition is given via an ``indirect method''. The main
result in \cite{Steen} shows that the data $(H^\cdot (\tilde
\fX^*,\C),L_\cdot, F^\cdot)$ determine a mixed $\Q$-Hodge
structure on $H^\cdot(\tilde\fX^*,\Q)$.

\smallskip

\no This result is obtained by studying  ``explicitly'' the
Picard--Lefschetz filtration $L_\cdot$ on a resolution of the
complex of sheaves $\Omega^\cdot_{\fX/\D}(\log Y)\
\otimes_{{\mathcal O}_{\fX}}{\mathcal O}_{Y^{\rm red}}$ of the
form
\begin{equation}\label{nearbyresol}
A^{s,k} := \Omega_{\fX}^{s+k+1}(\log Y)/W_k
\Omega_{\fX}^{s+k+1}(\log Y), \qquad s,k\in\Z_{\geq 0},
\end{equation}
where the weight filtration $W_\cdot$ on $\Omega_{\fX}^\cdot(\log Y)$ is
defined as
\[
W_k \Omega^m_{\fX}(\log Y) := \Omega^k_{\fX}(\log Y) \wedge
\Omega^{m-k}_{\fX},
\]
and the differentials
\begin{equation}\label{dApq}
d': A^{s,k} \to A^{s+1,k},\qquad d'': A^{s,k} \to A^{s,k+1}.
\end{equation}
on \eqref{nearbyresol} are  the usual differential $d'$ on
$\Omega^\cdot_{\fX}(\log Y)$ and
\begin{equation}\label{differ}
d''(\alpha) = (-1)^s\alpha\wedge\theta,\quad\text{for}\quad \theta
= f^*\left(\frac{dt}{t}\right).
\end{equation}
Notice that $\theta$ can be seen as an element of
\begin{equation}\label{H1Q1}
H^1(\fX^*,\Q)(1) = 2\pi\sqrt{-1}~H^1(\fX^*,\Q),
\end{equation}
because the form $dt/t$ on $\D^*$ has period $2\pi\sqrt{-1}$, so
that
\begin{equation}\label{H1ZQ}
(2\pi\sqrt{-1})^{-1} dt/t \in H^1(\D^*,\Z) \subset H^1(\D^*,\Q).
\end{equation}
Wedging with $\theta$ provides an injective map
\[
\wedge\theta: \Omega^m_{\fX/\D}(\log Y) \hookrightarrow
\Omega^{m+1}_{\fX}(\log Y)
\]
and an induced morphism of complexes of sheaves $ \phi:
\Omega^\cdot_{\fX/\D}(\log Y)\otimes_{{\mathcal O}_{\fX}}{\mathcal
O}_{Y^{\text{red}}} \to A^\cdot$, which defines a resolution of
the unipotent factor of the complex of nearby cycles, with
$(A^\cdot, \delta=d'+d'')$ the total complex of
\eqref{nearbyresol}. The endomorphism $\nu^{s,k}: A^{s,k} \to
A^{s-1,k+1}$ given by the natural projection on forms is
non-trivial because of the presence of the cutoff by the weight
filtration $W_k$ on $\Omega^m_{\fX}(\log Y)$ and it plays a
central role in this theory as it describes the local monodromy
map on the resolution $A^\cdot$.

\smallskip

\no It can be shown (see \cite{GNA})
that the map induced in hypercohomology by
\begin{equation}\label{tildenu}
\tilde\nu = (-1)^s\nu^{s,k}: A^{s,k} \to A^{s-1,k+1},
\end{equation}
satisfying $\tilde\nu \delta + \delta\tilde\nu = 0$, is the residue
\eqref{Nresidue} of the Gauss--Manin connection,
\[
N: H^m(\tilde \fX^*,\C) \to H^m(\tilde \fX^*,\C).
\]
More precisely, $N$ is obtained as the connecting homomorphism in
the long exact sequence of hypercohomology associated to the exact
sequence of complexes of sheaves on $Y$
\[
0 \to A^\cdot[-1] \stackrel{\epsilon}{\to}
\text{Cone}^\cdot(\tilde\nu) \stackrel{\eta}{\to} A^\cdot \to 0.
\]

\no The complex $(\text{Cone}^\cdot(\tilde\nu),D)$, $D = (\delta +
\tilde\nu,\delta)$ is quasi-isomorphic to $\Omega^\cdot_{\fX}(\log
Y)\otimes_{{\mathcal O}_{\fX}}{\mathcal O}_{Y^{\text{red}}}$. The
map $\tilde\nu$ measures the difference between differentiation in
$\text{Cone}^\cdot(\tilde\nu)$ and in $A^\cdot$ which appears when
one considers the section of $\eta$: $A^{s,k} \to
\text{Cone}^{s,k}(\tilde\nu) = A^{s-1,k}\oplus A^{s,k}$.

\smallskip
\subsection{Nearby cycles at arithmetic infinity}

\smallskip

\no The main idea that motivates the definition of the complex
\eqref{Kcomplex} is to ``transfer'' these results to the
archimedean setting, where one deals with a smooth, projective
algebraic variety $X$ over $\C$ or $\R$, interpreted as the
generic fiber of a degeneration ``around'' infinity. The main
aspects of the above construction that one wishes to retain are
the fact that the cutoff by $W_k$ on $\Omega^m_{\fX}(\log Y)$ is
introduced because the complex of the nearby cycles is the
restriction to $Y$ of the complex $\Omega^m_{\fX/\D}(\log Y)$. It
is the presence of this cutoff that makes the morphism $\nu^{s,k}$
non-trivial on $A^{s,k}$. Moreover, another essential observation
is that translates of the weight filtration
\begin{equation}\label{monof}
L_\ell A^{s,k} = W_{2k+\ell+1}\Omega_{\fX}^{s+k+1}(\log Y)/W_k
\Omega_{\fX}^{s+k+1}(\log Y);\qquad \ell\in\Z
\end{equation}
and the corresponding graded spaces
\begin{equation}\label{graded}
gr^L_\ell A^{s,k} =\left\{ \begin{array}{ll}
gr^W_{2k+\ell+1}\Omega_{\fX}^{s+k+1}(\log Y) & \ell+k \ge 0 \\[2mm]
0 & \ell+k<0 \end{array}\right.
\end{equation}
describe the strata of the special fiber $Y$, through the
Poincar\'e residue map
\begin{equation}\label{Pres}
{\rm  Res}: W_k \Omega^m_{\fX}(\log Y) \to
(a_k)_*\Omega^{m-k}_{\tilde Y^{(k)}},
\end{equation}
where
$$ \tilde Y^{(k)} := \displaystyle{\coprod_{1\le
i_1<\ldots<i_k\le M}}Y_{i_1}\cap\ldots\cap Y_{i_k} $$ is the
$k$-th stratum of $Y = Y_1\cup\ldots\cup Y_M$ and $(a_k)_*: \tilde
Y^{(k)} \to \fX$ is the canonical projection. The Poincar\'e
residue \eqref{Pres} is given by
\begin{equation}\label{Prescoord}
{\rm  Res} \left(\sum_{1\le i_1<\ldots<i_k\le K}\omega_{i_1\ldots
i_k}\frac{dz_{i_1}}{z_{i_1}}\wedge\ldots\wedge\frac{dz_{i_k}}{z_{i_k}}\right)
= \sum_{1\le i_1<\ldots<i_k\le K} {\rm res}(\omega_{i_1\ldots
i_k}),
\end{equation}
where $z_1\cdots z_K = 0$ is the local description of $Y^{\rm
red}$ (the closed subset $Y$ of $\fX$ with its reduced scheme
structure), $\omega_{i_1\ldots i_k}$ is a section of
$\Omega^{m-k}_{\fX}$, and ${\rm res}: \Omega_{\fX}^{m-k} \to
(a_k)_*\Omega^{m-k}_{\tilde Y^{(k)}}$ is the restriction to the
stratum $\tilde Y^{(k)}$.

\smallskip

\no This means that it is sufficient to provide a version at arithmetic
infinity of the weight filtration $W_\cdot$ and the Picard--Lefschetz
filtration $L_\cdot$, as these are sufficient to characterize the
geometry of the singular fiber $Y = f^{-1}(0)$, which is strictly
related to the behavior of the monodromy map.

\smallskip

\no Notice that the period $2\pi\sqrt{-1}$ of the form $dt/t$ on
$\D^*$ (\cf \eqref{differ}, \eqref{H1Q1}, \eqref{H1ZQ})
corresponds to a Tate twist on the (rational) cohomology of the
generic fiber $\tilde \fX^*$. It is important to stress the fact
that this `detects' the presence of the singular fiber through an
operation that does not involve $Y$ explicitly, and therefore can
be transported at arithmetic infinity (where the description of
the fiber `over infinity' is still mysterious) on a complex of
real differential forms. Moreover, the fact that the cutoff by
$W_k$ on $\Omega^m_{\fX}(\log Y)$ implies the non-triviality of
the monodromy map $\nu^{s,k}$ suggests a definition of the
monodromy operator $N$ ``at infinity'' in terms of an analogous
weight filtration on a complex of Tate--twisted real differential
forms.

\smallskip

\no By {\em weight} of a real $m$-form $$\alpha \in
\displaystyle{\bigoplus_{p+q=m}}(\Omega^{p,q}_X +
\Omega^{q,p}_X)_\R $$ on $X$ we mean the non-negative integer
$|p-q|$. At arithmetic infinity, the operations of taking the
residues and hence considering holomorphic differential forms on
the strata of $Y$ can be rephrased in a form suitable to be
included in the definition of the complex $K^\cdot$ of
\eqref{Kcomplex}.

\smallskip

\no In fact, one should interpret the archimedean complex
\eqref{Kcomplex} as a ``filtered copy'' of $A^{s,k}$, with an
additional condition characterizing the graded pieces $gr^L_\cdot
A^\cdot$. In fact, in the case of the complex of the nearby
cycles, there is a graded isomorphism
\begin{equation}\label{grac}
gr^L_\ell A^\cdot \simeq
\bigoplus_{k\ge\text{max}(0,-\ell)}(a_{2k+\ell+1})_*\Omega^{\cdot+1}_{\tilde
Y^{(2k+\ell+1)}}[-\ell-2k-1]
\end{equation}
under the condition $\ell+k\ge 0$. Furthermore, in $gr_\ell^L
A^\cdot$ the second differential is trivial, $d'' = 0$. The
induced {\em weights spectral sequence}
\begin{equation}\label{ss}
E_1^{i,m-i} = \bigoplus_{k\ge\text{max}(0,i)}H^{m+i-2k}(\tilde
Y^{(2k-i+1)},\Q)(i-k) \Rightarrow H^m(\tilde \fX^*,\Q)
\end{equation}
degenerates at the $E_2$ term. The $E_1^{i,m-i}$ term is a pure
Hodge structure of weight $m-i$. A major result in the theory
shows that the filtration induced on the abutment coincides with
the Picard-Lefschetz filtration.

\smallskip

\no At arithmetic infinity, in terms of the complex \eqref{Kcomplex},
for
$$ K^{i,j}=\bigoplus_{k\geq \max \{ 0, i \}} K^{i,j,k}, $$
the terms
\begin{equation}\label{E1infty}
K^{i,m-n,k}= \begin{cases} {\displaystyle
\bigoplus_{\begin{subarray}{l} p+q=m\\|p-q|\le 2k-i\end{subarray}}
(\Omega^{p,q}_X\oplus \Omega^{q,p}_X)_\R (\frac{m-i}{2})}
&\text{if $m-i \equiv~0(2),~k\ge\text{max}(0,i)$} \\
0 &\text{otherwise}
\end{cases}
\end{equation}
give the archimedean analog, at the level of the real differential
forms, of the $E_1^{i,m-i}$-term of the spectral sequence
\eqref{ss}. At arithmetic infinity the weight is $i-m$.

\smallskip

\no These analogies suggest a geometric interpretation of the
indices involved in the definition of \eqref{Kcomplex}. The first
index is associated to the $\ell$-th piece of an ``archimedean
monodromy filtration'' $L_{\ell}$ on
$\oplus_{p+q=m}(\Omega^{p,q}_X\oplus \Omega^{q,p}_X)_\R$, with
$\ell=-i$. This explains our previous comment that the archimedean
complex should be thought of as a filtered copy of $A^{s,k}$. The
cutoffs $|p-q|\le 2k-i$ and $k\ge \max\{ 0, i \}$ correspond to
considering the filtered piece $L_\ell
A^{s,k}=W_{2k+\ell+1}A^{s,k}$ in Steenbrink's theory with the
cutoff by $W_k$ on $\Omega_X^{s+k+1}(\log Y)$, justified by
restricting (relative) differential forms to the special fiber. In
this identification the third index $k$ of \eqref{Kcomplex} plays
the role of the index $k$ of the anti-holomorphic forms in the
double complex \eqref{nearbyresol}. The archimedean theory is a
weighted ``even theory'' since $i-m = 2r$. The second index
$j=m-n$ detects the total degree $m$ of the differential forms.

\smallskip

\no Finally, we remark that there is a fundamental difference in
the definition of the differentials \eqref{d1K} \eqref{d2K} in the
complex \eqref{Kcomplex} at infinity and their geometric analogs
\eqref{dApq}. In fact, while the differential $d'$ is similar in
both theories, in the geometric case, the action of $d''$ by the
wedging with the form $\theta = f^*(dt/t)$ involves a Tate twist
on the rational version of the complex $A^{s,k}$, the differential
$d''$ in the archimedean case does not involve any twist.

\section{Monodromy and the renormalization group}

\smallskip

\no In the construction at arithmetic infinity we obtain an analog of the
formula \eqref{Nresidue} for the logarithm of the monodromy $N$ as the
residue of a connection, in terms of a Birkhoff decomposition of loops
and a Riemann--Hilbert problem analogous to those underlying the
theory of renormalization in QFT, as developed by Connes and Kreimer
(\cite{CK2} \cite{CK3}, \cf also \cite{CoM}).

\smallskip

\no In our case, the Birkhoff decomposition will take place in the group
$G$ of automorphisms of the complex $(\T^\cdot_\C,\delta)$,
where $\T^\cdot_\C$ is the complexification of $\T^\cdot$. The
Birkhoff decomposition will determine a one parameter family of
principal $G$-bundles $\sP_\mu$
on $\bP^1(\C)$, with trivializations
\begin{equation}\label{birkhoff}
  \phi_\mu(z)=   \phi_\mu^-(z)^{-1}\, \phi^+_\mu(z), \ \ \  z\in
\partial\Delta \subset \bP^1(\C), \ \  \mu\in \C^*,
\end{equation}
where $\phi^+_\mu$ is a holomorphic function on a disk $\Delta$ around
$z=0$ and $\phi_\mu^-$ is holomorphic on $\bP^1(\C)\smallsetminus
\Delta$, normalized by $\phi^-_\mu(\infty)=1$. The parameter $\mu$ is
related to a scaling action by $\R^*_+$, by $\mu \mapsto \lambda
\mu$. We shall construct the data \eqref{birkhoff} in such a way that
the negative part of the Birkhoff decomposition $\phi^-_\mu= \phi^-$
is in fact independent of $\mu$, as in the theory of renormalization.
As part of the data, one also considers a one
parameter group of automorphisms $\theta$, and the corresponding
infinitesimal generator $\Upsilon =\frac{d}{dt}\theta_t |_{t=0}$, so
that
\begin{equation}\label{thetaphi1}
\phi_{\lambda \mu}(\epsilon) = \theta_{t\epsilon}\,
\phi_\mu(\epsilon), \ \ \ \  \forall \lambda=e^t \in \R^*_+, \ \  \epsilon\in
\partial\Delta.
\end{equation}
The corresponding {\em renormalization group} is the one parameter
semigroup
\begin{equation}\label{renormG}
\rho(\lambda)= \lim_{\epsilon \to 0 }\;\phi^-(\epsilon) \,
\theta_{t\epsilon} (\phi^- (\epsilon)^{-1}),  \ \ \ \  \forall \lambda
=e^t\in \R^*_+.
\end{equation}
Following \cite{CK2}, one can write $\phi^-(z)$ in the form
\begin{equation}\label{phiseries}
\phi^-(z)^{-1} = 1 + \sum_{k=1}^\infty  \frac{d_k}{z^k},
\end{equation}
with coefficients
\begin{equation}\label{dkcoeff}
d_k = \int_{s_1\ge \cdots \ge s_k \ge 0} \theta_{-s_1}(\beta) \cdots
\theta_{-s_k}(\beta) \, ds_1 \cdots ds_k.
\end{equation}
Here $\beta$ is the beta-function of renormalization, related to the
residue at zero of $\phi$ by
\begin{equation}\label{beta}
\beta = \Upsilon\, {\rm Res}\, \phi,
\end{equation}
where $\Upsilon$ is the generator of $\theta_t$ and the residue is
defined as
\begin{equation}\label{Resphi}
{\rm Res}\, \phi = \frac{d}{d z}\left( \phi^- (1/z)^{-1}
 \right)|_{z=0}.
\end{equation}

\smallskip

\no In the construction at arithmetic infinity,
the grading operator $\Phi$ induces a time evolution on the complex
$(\T^\cdot_\C,\delta)$ given by the ``Frobenius flow'' \eqref{Fflow},
\begin{equation}\label{Ft}
F_t = e^{t\Phi}, \ \ \  t\in \R,
\end{equation}
and we denote by $\theta_t$ the induced time evolution on
${\rm End}(\T^\cdot_\C,\delta)$,
\begin{equation}\label{thetat}
\theta_t(a) = e^{-t\Phi} \, a \, e^{t\Phi},
\end{equation}
with the corresponding $\Upsilon$ given by
\begin{equation}\label{YPhi}
\Upsilon(a) =\frac{d}{dt} \theta_t(a) |_{t=0} = [a,\Phi].
\end{equation}

\smallskip

\no The analogy with the complex of nearby cycles of a geometric
degeneration over a disk suggests to make the following prescription
for the residue of $\phi$:
\begin{equation}\label{Nresphi}
{\rm Res}\, \phi = N.
\end{equation}
As in the case of Connes--Kreimer, the residue uniquely determines
$\phi^-$ via \eqref{dkcoeff} and \eqref{phiseries}. We obtain the
following result.

\begin{The}\label{phiexpN}
There is a unique holomorphic map $\phi^-: \bP^1(\C)\smallsetminus \{
0 \} \to \Aut(\T^\cdot_\C)$ satisfying \eqref{phiseries}, with coefficients
\eqref{dkcoeff} and residue \eqref{Nresphi}, and it is of the form
\begin{equation}\label{expN}
\phi^-(z) = \exp(- N/z ).
\end{equation}
\end{The}

\proof First notice that the time evolution \eqref{thetat} satisfies
\begin{equation}\label{thetatN}
\theta_t(N)= e^t \,\, N.
\end{equation}
Moreover, by \eqref{YPhi} and \eqref{Nresphi}, we have
\begin{equation}\label{betaN}
\beta = [N,\Phi]=N
\end{equation}
Thus, we can write the coefficients $d_k$ of \eqref{dkcoeff} in the
form
$$ d_k = N^k \,\, \int_{s_1\ge \cdots \ge s_k \ge 0}
e^{-(s_1+\cdots+s_k)} \, ds_1 \cdots ds_k. $$
It is easy to see by induction that
$$
u_k(t):= \int_0^t \int_0^{s_1}\cdots \int_0^{s_{k-1}}
e^{-(s_1+\cdots+s_k)} \, ds_1 \cdots ds_k = \frac{(1-e^{-t})^k}{k!},
$$
satisfying the recursion $u_{k+1}^\prime (t) = e^{-t} u_k(t)$,
so that
$$ \int_{s_1\ge \cdots \ge s_k \ge 0}
e^{-(s_1+\cdots+s_k)} \, ds_1 \cdots ds_k = \frac{1}{k!}. $$
This implies that $d_k = N^k/k!$,
hence we obtain that the series \eqref{phiseries} is just
$$ \phi^-(z)^{-1} = \sum_{k=0}^\infty \frac{z^{-k}}{k!} N^k, $$
and $\phi^-(z)=\exp(-N/z)$.
\endproof

\smallskip

\no Correspondingly, we see that the renormalization group \eqref{renormG}
is given by
\begin{equation}\label{renN}
\rho(\lambda) = \lambda^N = \exp(t N), \ \ \ \  \forall \lambda=e^t\in
\R^*_+,
\end{equation}
since we have
\begin{equation}\label{thetatepphi}
\theta_{t\epsilon}(\phi^-(\epsilon)) = \exp\left(
-\frac{\lambda^\epsilon}{\epsilon}\, N \right) \ \ \ \  \forall
\lambda=e^t\in \R^*_+.
\end{equation}

\smallskip

\no We now show that the other term of the Birkhoff decomposition
\eqref{birkhoff} is determined by the requirement \eqref{thetaphi1} of
compatibility between the scaling and the time evolution.

\begin{The}\label{birkhoff0}
Consider the holomorphic map $\phi^+_\mu: \bP^1(\C)\smallsetminus \{
\infty \} \to \Aut(\T^\cdot_\C)$ given by
\begin{equation}\label{phi0}
\phi^+_\mu(z)= \exp\left( \frac{\mu^z-1}{z}\, N \right).
\end{equation}
The loop $\phi_\mu(z) =\phi^-(z)^{-1}\phi^+_\mu(z)$ with
$\phi^-(z)=\exp(-N/z)$ and $\phi^+_\mu$ as in \eqref{phi0}
satisfies the relation \eqref{thetaphi1}.
\end{The}

\proof First notice that \eqref{phi0} is indeed holomorphic at $z=0$
with $\phi^+_\mu(0)= \exp(\mu N)$.
By \eqref{birkhoff} and \eqref{thetatepphi}, the
relation \eqref{thetaphi1} is equivalent to requiring that, for all
$\lambda=e^t\in
\R^*_+$, the function $\phi^+_{\mu}$ satisfies
\begin{equation}\label{evoleqphi}
\phi^+_{\lambda\mu}(\epsilon)= \exp\left(
\frac{\lambda^\epsilon-1}{\epsilon}\, N \right) \, \theta_{t\epsilon}(
\phi^+_{\mu}(\epsilon)).
\end{equation}
For $\phi_\mu^+$ as in \eqref{phi0} we have
$$ \theta_{t\epsilon}(\phi^+_{\mu}(\epsilon))= \exp\left(
\frac{\mu^\epsilon-1}{\epsilon}\, \lambda^\epsilon\, N \right) =
\exp\left( \frac{(\lambda\mu)^\epsilon-1}{\epsilon}\, N \right)
\exp\left(\frac{1-\lambda^\epsilon}{\epsilon}\, N \right), $$
so that \eqref{evoleqphi} is satisfied.
\endproof

\smallskip

\no Notice that, via the representation \eqref{SL2R2}, it is possible to
lift the Birkhoff decomposition \eqref{birkhoff} to a Birkhoff
decomposition of the form
\begin{equation}\label{birkSL2}
g_\mu(z) = g^-(z)^{-1} g^+_\mu(z),
\end{equation}
where, with the notation of \eqref{SL2Relts}, we have
$$ g^-(z) = u(1/z), \ \ \  g_\mu(z)= u(\mu^z /z). $$
The renormalization group \eqref{renN} then
becomes simply the horocycle flow
\begin{equation}\label{renhoro}
\rho(\lambda)= 
u(t)= \left(\begin{array}{cc} 1 & t \\ 0 & 1
\end{array} \right).
\end{equation}

\smallskip

\no There is a Riemann--Hilbert problem associated to the Birkhoff
decomposition considered in the theory of renormalization (\cf
\cite{CoM}). Namely, for $\gamma$ a generator of the fundamental
group $\pi_1(\Delta^*)=\Z$ of the punctured disk, consider a
complex linear representation $\pi: \Z \to G$. Under the
assumption that the eigenvalues of $\pi(\gamma)$ satisfy
$$ 0\leq {\rm Re}\, \frac{\lambda}{2\pi\sqrt{-1}} <1, $$
we can take the logarithm
\begin{equation}\label{repRH}
\frac{1}{2\pi\sqrt{-1}}\, \log \pi(\gamma).
\end{equation}
By the Riemann--Hilbert correspondence (\cf \cite{AnBo}), a
representation $\pi: \Z \to G$ determines a bundle with connection
$(\sE,\nabla)$, where the Fuchsian connection $\nabla$ has local gauge
potential on the disk $\Delta$ of the form
\begin{equation}\label{gaugeA}
 -\phi^+(z)^{-1} \frac{\log \pi(\gamma)\, dz}{z} \phi^+(z) +
\phi^+(z)^{-1}\, d\phi^+(z),
\end{equation}
with $\phi^+(z)$ the local trivialization of $\sE$ over $\Delta$.
The data $(\sE,\nabla)$ correspond to a linear differential system
$f'(z) = A(z) f(z)$, with $\nabla f = df - A(z) f dz$, for
sections $f\in \Gamma(\Delta^*, \sE)$.

\smallskip

\no In the case of the cohomological theory at arithmetic infinity, this
amounts to a vector bundle $\sE_\mu^\cdot$ over $\bP^1(\C)$, with
fiber $\T^\cdot_\C$, associated to the principal $G$-bundles $\sP_\mu$,
with transition function the loop $\phi_\mu(z)$ and local
trivializations given by the $\phi^-(z)$ and $\phi_\mu^+(z)$. We use
the representation specified by
\begin{equation}\label{piT0}
\pi(\gamma):= \exp(-2\pi\sqrt{-1} \,\, N),
\end{equation}
which is the analog of \eqref{expNres}.
This determines a one parameter family $(\sE_\mu, \nabla_\mu)$ of
linear differential systems over the disk $\Delta$, of the form
\begin{equation}\label{Fuchsian}
\nabla_\mu : \sE_\mu^\cdot \to \sE_\mu^\cdot \otimes_{\sO_\Delta}
\Omega^\cdot_{\Delta}(\log 0),
\end{equation}
where the connection on the restriction of $\sE_\mu$ over $\Delta$ is
given by
\begin{equation}\label{connectionE}
\nabla_\mu = N\, \left( \frac{1}{z} + \frac{d}{dz} \frac{\mu^z-1}{z}
\right) \, dz.
\end{equation}
Using the induced representation \eqref{SL2R2} in cohomology, this
determines a corresponding linear differential system on the bundle of
hypercohomologies $\H^\cdot(\sE_\mu^\cdot)$ with fiber
$\H^\cdot(\T^\cdot_\C)$,
\begin{equation}\label{FuchsianH}
\nabla_\mu : \H^\cdot(\sE_\mu^\cdot) \to \H^\cdot(\sE_\mu^\cdot)
\otimes_{\sO_\Delta} \Omega^\cdot_{\Delta}(\log 0),
\end{equation}
which is the analog of the Gauss--Manin connection \eqref{GaussManin}
in the geometric case. The connections \eqref{connectionE} form
an isomonodromic family, with
\begin{equation}\label{isomonodromy}
{\rm Res}_{z=0} \nabla_\mu =N.
\end{equation}

\section{Rees sheaves at arithmetic infinity}

\smallskip

\no In the description \eqref{DenDet} (\cf \cite{Den}) of the
archimedean factor of the Hasse-Weil L-function of the ``motive''
$H^m(X)$, for $X$ a smooth projective variety over a number field,
the definition of the archimedean cohomology $H^m_{ar}(X)$ is
motivated by previous work of J.~M.~Fontaine and it is expressed
in terms of an additive functor $\mathbb D$ (derivation) from the
(abelian) category of pure Hodge structures over $\kappa = \C,\R$
to the additive category whose objects are free modules of finite
rank over the $\R$-algebra of polynomials in one variable endowed
with a $\R$-linear endomorphism and satisfying certain properties.
More precisely one sets
\[
H^m_{ar}(X) = \begin{cases} Fil^0(H^m_B(X,\C)\otimes_\C \C[z^{\pm
1}])^{c=id}  &\text{if $\kappa = \C$} \\
Fil^0(H^m_B(X,\C)\otimes_\C \C[z^{\pm 1}])^{c=id, F_\infty = id}
&\text{if $\kappa = \R$}.  \end{cases}
\]
Here, $Fil^q$ denotes the filtration on the tensor product
$H^m_B(X(\C),\C)\otimes_\C \C[z^{\pm 1}]$ obtained from the Hodge
filtration $F^\cdot$ on the Betti cohomology $H^m_B(X) :=
H^m_B(X(\C),\C)$ and the one on the ring of Laurent polynomials
$\C[z^{\pm 1}]$ given by $F^q\C[z^{\pm 1}] :=
z^{-q}\C[z^{-1}],~\forall q\in\Z$. By $c$ one denotes the
conjugate linear involution (complex conjugation) and $F_\infty$
is a $\C$-linear involution (the infinite Frobenius).

\smallskip

\no The expectation is that one should obtain a description of the
archimedean cohomology together with the linear ``Frobenius flow'' 
generated by $\Phi$
directly by some natural homological construction on a suitable
non-linear dynamical system. In this direction, a more geometric
construction of the archimedean cohomology was given in
\cite{Den2} (\cf~par.~3), where it was interpreted (for instance
when $\kappa =\C$) as the space of global sections of a real analytic
sheaf (Rees sheaf) $\zeta_{\C}^\omega(H^m(X),\gamma^\cdot)$ over $\R$. A
similar description holds when $\kappa = \R$. Here, as before, 
$\gamma^\cdot$
denotes the descending filtration $F^\cdot\cap\bar F^\cdot$ on
$H^m(X^{\rm{an}},\R)$ endowed with its real Hodge structure. The
locally-free sheaf $\zeta_{\C}^\omega(H^m(X),\gamma^\cdot)$ has
the remarkable description (\cf~\cite{Den2} Thm.~4.4)
\begin{align}\label{Rsheaf}
\zeta_{\C}^\omega(H^m(X),\gamma^\cdot) &\simeq
\Ker \left(\R^m\pi_*(\Omega^\cdot_{X^{\rm{an}}_\C\times\R/\R},sd)
\to (\R^m\pi_*\mathcal{DR}_{X/\C})/\mathcal T_{X/\C}\right)\\
H^m_{\rm{ar}}(X) &=
\Gamma(\R,\zeta_{\C}^\omega(H^m(X),\gamma^\cdot) \qquad
\text{if}\quad \kappa = \C \notag
\end{align}
where $\pi: X^{\rm{an}}_\C\times \R \to \R$ is the projection, $s$
is a standard coordinate on $\R$, $\mathcal A_\R$ denotes the
sheaf of real-analytic functions on the real analytic manifold
$\R$, $\mathcal{DR}_{X/\C}$ is the cokernel of the natural
inclusion of complexes of $\pi^{-1}\mathcal A_\R$-modules on
$X^{\rm{an}}_\C\times \R$
\[ \pi^{-1}\mathcal A_\R \hookrightarrow
(\Omega^\cdot_{X^{\rm{an}}_\C\times\R/\R},sd)
\]
\no and $\mathcal T_{X/\C}$ is the $\mathcal A_\R$-torsion in
$\R^m\pi_*\mathcal{DR}_{X/\C}$.
$\Omega^\cdot_{X^{\rm{an}}_\C\times\R/\R}$ is the complex of
$\C$-valued smooth relative differential forms on
$X^{\rm{an}}\times\R/\R$ which are holomorphic in the
$X^{\rm{an}}$-coordinates and real analytic in the $\R$-variable.
One considers the scaling flow on $\R$ given by the map $\phi^t_\C(s) =
se^{-t}$. It induces an action of this group on the relative
differential complex by means of: $\psi^t(\omega) =
e^{t\rm{deg}\omega}\cdot(\rm{id}\times \phi^t_\C)^*\omega$. This
action defines in turn a $\mathcal A_\R$-linear action on the
complex of higher direct image sheaves on $\R$. A similar result
holds when $\kappa = \R$: in this case one gets a $\mathcal
A_{\R^{\ge 0}}$-action. In the description of $H^m_{\rm{ar}}(X)$
given in \eqref{Rsheaf}, the Hodge theoretic notions required in
the definition of the archimedean cohomology have been replaced by
using suitably deformed complexes of sheaves of modules on $\R$
(on $\R^{\ge 0}$ when $\kappa = \R$). The deformed complex of
locally free sheaves of relative differentials
$(\Omega^\cdot_{X_\C\times\bA^1/\bA^1},zd)$, filtered by the Hodge
filtration $F^\cdot$, for $z$ coordinate on $\bA^1$, was firstly
studied by Simpson in \cite{Simpson}. He introduced the algebraic
version $\zeta_\C(H^m(X^{\rm{an}},\C),F^\cdot)$ of the real
analytic Rees sheaf and proved that
$\zeta_\C(H^m(X^{\rm{an}},\C),F^\cdot) \simeq
\R^m\pi_*(\Omega^\cdot_{X_\C\times\bA^1/\bA^1},zd)$. Following
this viewpoint, \eqref{Rsheaf} is the analogue of Simpson's
formula for the non-algebraic filtration $\gamma^\cdot$. A very
interesting fact (\cf~\cite{Den2}, \S 4) is that in the real
analytic setting, the higher direct image sheaves fit into short
exact sequences of coherent $\mathcal A_\R$-modules
\begin{equation}\label{ChD}
0\to \R^m\pi_*(\pi^{-1}\mathcal A_\R) \to
\R^m\pi_*(\Omega^\cdot_{X^{\rm{an}}_\C\times\R/\R},sd)
\stackrel{\alpha}{\to} \R^m\pi_*\mathcal{DR}_{X/\C}\to 0
\end{equation}

\no where $\R^m\pi_*(\pi^{-1}\mathcal A_\R) = H^m(X^{\rm{an}},\R)
\otimes \mathcal A_\R$. We like to think the sheaf
$\R^m\pi_*(\Omega^\cdot_{X^{\rm{an}}_\C\times\R/\R},sd)$ as the
archimedean real analytic analog of the relative analytic
hypercohomology $\R^mf_*\Omega^\cdot_{\mX/\Delta}(\log Y)$ over a
small disk $\Delta$ centered at the origin, whose algebraic
description has been recalled in \S~3.3. It turns out that
$\zeta_{\C}^\omega(H^m(X),\gamma^\cdot)$ is also canonically
isomorphic to a twisted dual of $\R^m\pi_*\mathcal{DR}_{X/\C}$.
More precisely, if $d,m\in\Z_{\ge 0}$ with $m+d = 2n$ ($n = \dim
X$), then there are isomorphisms of $\mathcal A_\R$-modules
(\cf~\cite{Den2} Thm.~4.2)
\begin{equation}\label{Rsheaf2}
\zeta_{\C}^\omega(H^m(X),\gamma^\cdot) \simeq
(2\pi\sqrt{-1})^{1-n}\mathcal Hom_{\mathcal
A_\R}(\R^d\pi_*\mathcal {DR}_{X/\C},\mathcal A_\R(-n)).
\end{equation}
\smallskip
\no The dualizing operation detects the $\gamma^\cdot$-filtration
from the Hodge filtration on $\R^d\pi_*\mathcal {DR}_{X/\C}$. From
\eqref{Rsheaf} and \eqref{Rsheaf2} one gets isomorphisms
respecting the $\mathcal A_\R$-module structures and the flow
\begin{equation}\label{dual}
\Ker\left(\R^m\pi_*(\Omega^\cdot_{X^{\rm{an}}_\C\times\R/\R},sd)
\to (\R^m\pi_*\mathcal{DR}_{X/\C})/\mathcal T_{X/\C}\right) \simeq
(2\pi\sqrt{-1})^{1-n}\mathcal Hom_{\mathcal
A_\R}(\R^d\pi_*\mathcal {DR}_{X/\C},\mathcal A_\R(-n)).
\end{equation}

\no This statement is the dynamical sheaf-theoretic analog of the
duality isomorphisms between the hypercohomology of the complexes
of vector spaces $\Ker(N)^\cdot$ and $\Coker(N)^\cdot$ of
\cite{KC} (Prop.~4.13) induced by powers of the `local monodromy
at arithmetic infinity', that is, of the duality $S$ of \eqref{involution}.

\no This way, one can reinterpret the
archimedean cohomology as the space of global sections on $\R$ (on
$\R^{\ge 0}$ when $\kappa = \R$) of the sheaf inverse image in
$\R^m\pi_*(\Omega^\cdot_{X^{\rm{an}}\times\R/\R},sd)$ of the
maximal $\mathcal A_\R$-submodule of
$\R^m\pi_*\mathcal{DR}_{X/\C}$ with support in $0\in\R$. This
statement is in accord with the classical description of the
inertia invariants as the kernel of the local monodromy map $N$,
viewed as the residue at zero of the Gauss-Manin
connection
$$\Ker ({\rm Res}_0\nabla: \R^m
f_*\Omega^\cdot_{\mX/\Delta}(\log Y)\otimes_{\mathcal O_\Delta}
k(0) \to \R^mf_*\Omega^\cdot_{\mX/\Delta}(\log Y)\otimes_{\mathcal
O_\Delta} k(0) )$$ \cf~\S~3.3 and \eqref{Fuchsian},
\eqref{connectionE}, \eqref{FuchsianH} in the archimedean case.

\no In the archimedean setting (\eg for $\kappa = \C$), the exact
sequences
\begin{equation}\label{arseq}
0 \to \zeta^\omega_\C(H^m(X),\gamma^\cdot) \to
\R^m\pi_*(\Omega^\cdot_{X^{\rm{an}}_\C\times\R/\R},sd) \to
(\R^m\pi_*\mathcal{DR}_{X/\C})/\mathcal T_{X/\C}\to 0
\end{equation}
\no are the sheaf theoretic analogs of the hypercohomology exact
sequences associated to the nearby cycles complex as defined by
Deligne in \cite{SGA7II} (\cf Exp.~XIII, \S~1.4: (1.4.2.2)). At
arithmetic infinity, the hypercohomology of the complex of
vanishing cycles is replaced  by the hypercohomology sheaf
$(\R^m\pi_*\mathcal{DR}_{X/\C})/\mathcal T_{X/\C}$, whereas the
archimedean analog of the variation map in the formalism of the
nearby cycles produces the duality isomorphisms \eqref{dual}.

\smallskip

\no The sequence \eqref{arseq} has also interesting analogies with
the exact sequence of $C_k$-modules ($C_k =$ id\`ele class group
of a global field $k$)
\begin{equation}\label{AlC}
0 \to L^2_\delta(X)_0 \stackrel{E}{\to} L^2_\delta(X) \to \mathcal
H \to 0
\end{equation}

\no studied by A.~Connes in \cite{AC}, \S III (33). It is
tempting to interpret the P\'olya-Hilbert space $\mathcal H$ as
the Hilbert space analogue of the cohomology
$(\R^\cdot\pi_*\mathcal{DR}_{X/\C})/\mathcal T_{X/\C}$. This
relations suggest a `singular behavior' of $\overline{{\rm
Spec}(\Z)}$ around infinity.
\smallskip

One of the interesting questions related to the archimedean
cohomology and the archimedean factor is that of writing the
logarithm $\log L_\kappa(H^m(X),s)$ of the regularized determinant
\eqref{DenDet} via a Lefschetz trace formula for the Frobenius
operator. To this purpose, a first necessary step seems that of
relating the complex $\T^\cdot$ to a complex of sheaves of
(relative) differential forms on $X^{\rm{an}}_\C\times \R$ and
eventually interpreting the Tate twists in the complex $K^\cdot$
(the degree in the variable $U$ in $\T^\cdot$) as ``Fourier
modes''. We hope to develop these topics in a future work.

\section{``Arithmetic'' spectral triples}

\smallskip

\no In this section we present a refined version of the 
proposed construction of an
``arithmetic spectral triple'' in \cite{CM}. This version has the
advantage that it holds at the level of differential forms and for
$X$ of any dimension. 
We first introduce natural subcomplexes and quotient complexes 
of $\T^\cdot$.

\smallskip

\subsection{Inertia invariants and coinvariants}

\smallskip

\begin{figure}
\begin{center}
\epsfig{file=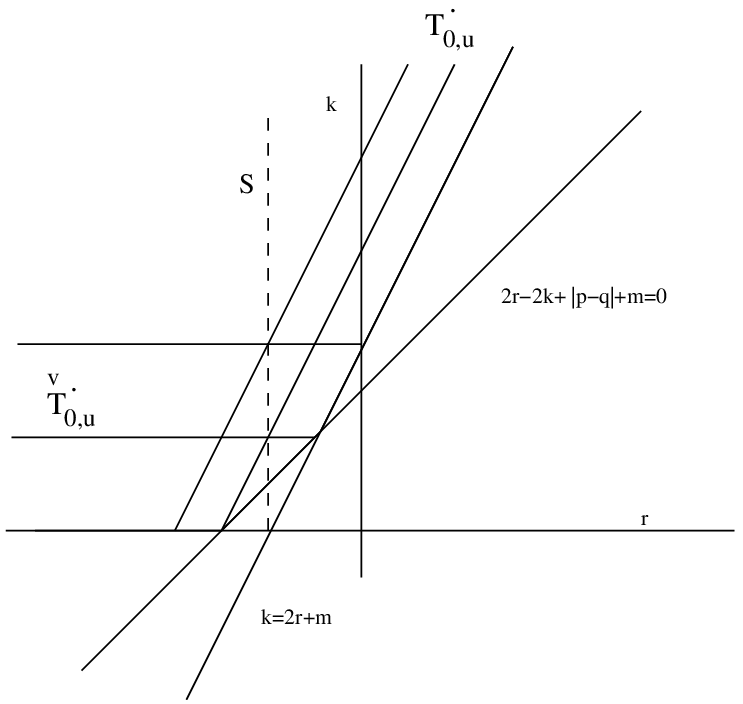} \caption{The complexes
$\check\T_{0,u}^\cdot$ and $\T^\cdot_{0,u}$ and the duality $S$
\label{TLFig}}
\end{center}
\end{figure}

\no We consider certain complexes of vector spaces related to the
action of the endomorphism $N$ on $(\T^\cdot,\delta)$. We
introduce the notation $\T^\cdot_\ell\subset \T^\cdot$ for the
$\Z$-graded linear subspace obtained as follows. For fixed $(p,q)$
with $p+q=m$, let $(\T^{m,*}_{p,q})_\ell \subset \T^{m,*}_{p,q}$
be the $2\Z$-graded real vector space spanned by elements of the
form $\alpha \otimes U^r \otimes \hbar^k$, with $\alpha \in
(\Omega^{p,q}_X\oplus \Omega^{q,p}_X)_\R$ and $(r,k)\in
\Lambda_{p,q}$ lying on the line $k=2r+m+\ell$ (\cf Figure
\ref{TLFig}). We let $\T^\cdot_\ell=\oplus_{p,q}
(\T^{m,*}_{p,q})_\ell$.

\smallskip

\no Each $\T^\cdot_\ell$ is a subcomplex with respect to the
differential $d'=\hbar d$, while the second differential satisfies
$d'' : \T^\cdot_\ell \to \T^{\cdot\, +1}_{\ell-1}$.

\smallskip

\no Similarly, we denote by $\check\T^\cdot_\ell$ the linear
subspace of $\T^\cdot$, which is the direct sum of the subspaces
$(\check\T^{m,*}_{p,q})_\ell\subset \T^{m,*}_{p,q}$ spanned by
elements $\alpha \otimes U^r \otimes \hbar^k$, with $\alpha \in
(\Omega^{p,q}_X\oplus \Omega^{q,p}_X)_\R$ and with $(r,k)\in
\Lambda_{p,q}$ lying on the  horizontal line $k=\ell$ (\cf Figure
\ref{TLFig}). Each $\check\T^\cdot_\ell$ is a subcomplex with
respect to the differential $d''$, while $d':\check\T^\cdot_\ell
\to \check\T^{\cdot+1}_{\ell+1}$.

\smallskip

\no In particular, for $\ell=0$ and $k=2r+m\geq 0$, we obtain the
subcomplex of ``inertia invariants''
\[
\T^\cdot_0 := \Ker(N)^\cdot. \]

\no This complex agrees with the complex that computes the
$\H^\cdot(\T^\cdot,\delta)^{N=0}$, for all $k=2r+m> 0$, so that
the map
$$ \H^*(\Ker(N)^\cdot,d) \twoheadrightarrow \H^*(\T^\cdot,d)^{N=0} $$
is almost always a bijection.

\smallskip

\no Similarly, for $u\in \N$, one can consider subcomplexes
$\T_{0,u}^\cdot\subset \T^\cdot$
\begin{equation}\label{TL}
\T_{0,u}^\cdot :=\oplus_{\ell=0}^u \T^\cdot_\ell=
\Ker(N^{u+1})^\cdot.
\end{equation}
These satisfy $\T^\cdot = \varinjlim_u \T^\cdot_{0,u}$. We also
consider the quotient complexes
\begin{equation}\label{cTL}
\check\T_{0,u}^\cdot : =\oplus_{\ell=0}^u \check\T^\cdot_\ell =
\Coker(N^{u+1})^\cdot,
\end{equation}
where we denote by $\check\delta$ the induced differential on
$\check\T^\cdot_{0,u}$. We obtain this way $\T^\cdot =
\varprojlim_u \check\T^\cdot_{0,u}$. We call $\Coker(N)^\cdot$ the
complex of ``inertia coinvariant'' and we refer to the
$\Ker(N^{u+1})^\cdot$ and $\Coker(N^{u+1})^\cdot$ as higher
inertia co/invariants.

\bigskip

\subsection{Spectral triple}

\smallskip

\no We consider the complex
\begin{equation}\label{T-truncated}
\T^\cdot_u = \T^\cdot_{0,u} \oplus \check\T^\cdot_{0,u}[+1],
\end{equation}
with induced differential $\delta_u=\delta\oplus \check\delta$,
where $(\T^\cdot_{0,u},\delta)$ and
$(\check\T^\cdot_{0,u},\check\delta)$ are the complexes of higher
invariants and coinvariants of \eqref{TL} and \eqref{cTL}.

\no We still denote by $\Phi$ the linear operator on $\T^\cdot_u$
which agrees with the operator $\Phi=-U\partial_U$ on the subspaces
$0\oplus \check\T^\cdot_{0,u}$ and $\T^\cdot_{0,u} \oplus 0$,
identified with the corresponding linear subspaces of $\T^\cdot$.

\no On the compact K\"ahler manifold $X$ consider the operator $d+d^*$ on
real forms $\Omega^\cdot_{X,\R}$. For the next result we need to
assume the following conditions on the spectrum of $d+d^*$:
\begin{equation}\label{SP-cond}
\begin{array}{ll}
\# \{ (r,\lambda)\in \Z\times \Sp(d+d^*): \, t=r+\lambda \} < \infty &
\forall t\in \R \\[3mm]
\{ r+\lambda\,: \,\, r\in \Z,\, \lambda \in \Sp(d+d^*) \}\subset \R &
\text{is discrete}.
\end{array}
\end{equation}

\begin{The}\label{arith3}
Let $U(\cg)$ be the universal enveloping algebra of the Lie algebra
$\cg=\gsl(2,\R)$. Let $\sH^\cdot$ denote the completion of
$\T^\cdot_u$ with respect to the inner product induced by
\eqref{innprod}. Let $\sD$ be the linear operator $\sD=\Phi +
\delta_u +\delta_u^*$. Let $\cA=\sC^\infty(X,\R)\otimes
U(\cg)$. If $X$ has the property \eqref{SP-cond}, then
the data $(\cA,\sH^\cdot, \sD)$
satisfy the properties:
\begin{itemize}
\item The representation \eqref{SL2L2} determines an action of the
algebra $\cA$ by bounded operators on the Hilbert space $\sH^\cdot$.
\item The commutators $[\sD,\sigma^L(a)]$, for $a\in \cA$, are
bounded operators on $\sH^\cdot$.
\item The operator $\sD$ is a densely defined unbounded
self-adjoint operator on $\sH^\cdot$, such that
$(1+\sD^2)^{-1}$ is a compact operator.
\end{itemize}
\end{The}

\proof The representation $\sigma^L$ of $\SL(2,\R)$ on $\T^\cdot$
defined in \eqref{SL2L2} preserves the subspaces $\T^\cdot_{0,u}$ and
$\check\T^\cdot_{0,u}$. In fact, the operator $\hL$ changes $m\mapsto
m+2$, $r\mapsto r-1$ and $k\mapsto k$, so that the constraint $k\leq
u$ defining $\check\T^\cdot_{0,u}$ and the constraint $k=2r+m+\ell$
with $0\leq \ell\leq u$ defining $\T^\cdot_{0,u}$ are preserved by the
action of $\hL$. Similarly, the involution $\tilde S$ changes
$m\mapsto 2n-m$, $r\mapsto r-(n-m)$ and $k\mapsto k$, so that again
both constraints $k\leq u$ and $k=2r+m+\ell$, for
$0\leq \ell\leq u$, are preserved. Thus, we can consider on
$\T^\cdot_u$ the corresponding derived representation $d\sigma^L$ of
the Lie algebra $\cg=sl(2,\R)$,
$$ d\sigma^L(v)=\frac{d}{ds} \sigma^L(\exp(sv))|_{s=0}. $$
Let $\{ v_{\pm},v_0 \}$ be the basis of $\cg$ with $[v_0,v_+]=2v_+$,
$[v_0,v_-]=-2v_-$, and $[v_+,v_-]=v_0$. We have $d\sigma^L(v_+)=\hL$,
while $d\sigma^L(v_0)$ is the linear operator that multiplies elements
$\alpha\otimes U^r\otimes \hbar^k$ with $\alpha\in \Omega^m_{X,\R}$
by $-n+m$. Thus, we obtain an action of the algebra $U(\cg)$ on
$\sH^\cdot$ by bounded linear operators.

\no We have
$[\sigma^L(\gamma),\delta_u]=[\sigma^L(\gamma),\delta_u^*]=0$ hence
$[\sD,d\sigma^L(v)]=[\Phi,d\sigma^L(v)]$. Using $[\Phi,\hL]=\hL$ and
$[\Phi,\tilde S]=(-n+m) \tilde S$, we obtain that $[\sD,d\sigma^L(v)]$
is a bounded operator for all $v\in U(\cg)$. The algebra of real
valued smooth functions $\sC^\infty(X,\R)$ acts on $\T^\cdot$ and on
$\T^\cdot_u$ by the usual action on real forms
$\Omega^\cdot_{X,\R}$. This action commutes with $\Phi$ so that
$[\sD,f]=[\delta_u+\delta_u^*,f]$, for all $f\in \sC^\infty(X,\R)$,
and we can estimate $\| \, [\sD,f] \, \|\leq C \sup |df|$ for some
$C>0$.

\no For $\sD_u=\delta_u+\delta_u^*$, we have $\sD=\Phi+\sD_u$, with
$[\Phi,\sD_u]=0$. The operator $\Phi$ on $\T^\cdot_u$ has spectrum
$\Z$. The eigenspace $E_r$ with eigenvalue $r\in \Z$ is the span
of the elements
$$ (\alpha_1\otimes U^r\otimes \hbar^{k_1}, \alpha_2\otimes U^r\otimes
\hbar^{k_2}) \in \T^\cdot_{0,u}\oplus \check\T^\cdot_{0,u}, $$
namely, elements with $\alpha_i\in \Omega^{m_i}_{X,\R}$ and $k_i$
satisfying $0\leq k_1-2r-m_1 \leq u$ and $0\leq k_2 \leq u$.
The operator $\sD_u$ restricted to the eigenspace $E_r$ has discrete
spectrum. The multiplicity $m_\lambda$ of an eigenvalue $\lambda$ of
$\sD_u|_{E_r}$ is bounded by $2u\, n_\lambda$, where $n_\lambda$ is the
multiplicity of $\lambda$ as an eigenvalue of the operator $d+d^*$ on
real differential forms $\Omega^\cdot_{X,\R}$. Condition
\eqref{SP-cond} then implies that $\sD$ has discrete spectrum with
finite multiplicities and that $(1+\sD^2)^{-1}$ is compact.

\endproof

\no When passing to cohomology, we obtain an induced structure of the
following form.

\begin{Cor}\label{arith3Cor}
Let now $\sH^\cdot$ denote the Hilbert completion of the cohomology
$\H^\cdot(\T^\cdot_u,\delta_u)$ and let $\sD=\Phi$, the
operator induced in cohomology. For $\cA=U(\cg)$, the triple
$(\cA,\sH^\cdot, \sD)$ satisfies:
\begin{itemize}
\item The representation \eqref{SL2L2} induces an action of
algebra $\cA$ by bounded operators on $\sH^\cdot$.
\item The commutators $[\sD,\sigma^L(a)]$, for $a\in \cA$, are
bounded operators on $\sH^\cdot$.
\item The operator $\sD$ is a densely defined unbounded
self-adjoint operator on $\sH^\cdot$,
such that $(1+\sD^2)^{-1}$ is a compact operator.
\end{itemize}
\end{Cor}

\proof The statement follows as in the case of Theorem
\ref{arith3}. Notice that now we have $\Sp(\Phi)=\Z$ with
multiplicites
$$ m_r = \dim \frac{\Ker( d': \T^{\cdot,2r}_{0,u}
\to \T^{\cdot+1,2r}_{0,u})}{{\rm Im}(d':\T^{\cdot-1,2r}_{0,u}
\to \T^{\cdot,2r}_{0,u})}
+ \dim \frac{\Ker( d'': \check\T^{\cdot+1,2r}_{0,u}
\to \check\T^{\cdot+2,2r}_{0,u})}{{\rm Im}(d'':\check\T^{\cdot,2r}_{0,u}
\to \check\T^{\cdot+1,2r}_{0,u})}. $$
This cohomological result no longer depends on the property \eqref{SP-cond}.

\endproof

\no Recall that $(\cA,\sH^\cdot, \sD)$ is a spectral triple in the 
sense of Connes (\cf \cite{Co3}) if the three properties 
listed in the statement of Theorem \ref{arith3} and Corollary 
\ref{arith3Cor} hold and the algebra $\cA$ is a dense involutive 
subalgebra of a $C^*$-algebra.

\medskip

\no In our case, the adjoints of elements in $\cA$ with 
respect to the inner product on $\sH^\cdot$ are again contained in $\cA$. 
In fact, one can see that the adjoint of the Lefschetz $\hL$ is given by
\begin{equation}\label{Lef-adj}
\hL^* = (\cdot \rfloor \omega )\, U,
\end{equation}
where $\rfloor$ is the interior product and $\omega$ is the K\"ahler
form, and we obtain \eqref{Lef-adj} from
$$ \sigma^L(w)^{-1} \hL \sigma^L(w) = *((*\cdot)\wedge \omega)\, U. $$

\smallskip

\no Moreover, a choice of the K\"ahler form (of the K\"ahler class
in the cohomological case) determines a corresponding representation
of the involutive algebra $\cA$ on the Hilbert space $\sH^\cdot$.
The choice of the K\"ahler class ranges over the K\"ahler cone
\begin{equation}\label{Kcone}
\sK = \{ c\in H^{1,1}(X): \,\, \int_M c^k >0 \}
\end{equation}
for all $M\subset X$ complex submanifolds of dimension $1\leq k\leq
X$. Thus, in the case of Corollary \ref{arith3Cor} for
instance we can consider a norm
\begin{equation}\label{norm}
\| a \|:= \sup_{c\in \overline{\sK}:\, \|c\|=1} \| \sigma^L_c(a) \|,
\end{equation}
where $\overline{\sK}$ is the nef cone and $\sigma^L_c$ is the
representation of $\cA$ determined by the choice of the class $c$.
Thus, the data $(\cA,\sH^\cdot, \sD)$ of Theorem \ref{arith3}
and Corollary \ref{arith3Cor} determine a spectral triple.

\smallskip

\no An interesting arithmetic aspect of spectral triples is that 
they have an associated family of zeta functions, the 
basic one being the zeta function of the Dirac operator, 
$$ \zeta_{\sD} (z) = \Tr(|\sD|^{-z}) = \sum_{\lambda} 
\Tr(\Pi(\lambda,|\sD|)) \lambda^{-z}, $$
where $\Pi(\lambda,|\sD|)$ denotes the orthogonal projection onto the
eigenspace $E(\lambda,|\sD|)$. More generally, one considers 
for $a\in \cA$ the corresponding zeta function
$$ \zeta_{a,\sD}(z) = \Tr (a |\sD|^{-z}) = \sum_{\lambda}
\Tr(a\, \Pi(\lambda,|\sD|)) \lambda^{-z}. $$
These provide a refined notion of dimension for noncommutative spaces, 
the dimension spectrum, which is the complement in $\C$ of the set where
all the $\zeta_{a,\sD}$ extend holomorphically.

\no One can also consider the associated two-variable zeta functions,
$$ \zeta_{a,\sD}(s,z) := \sum_{\lambda} \Tr(a\, \Pi(\lambda,|\sD|))
(s-\lambda)^{-z} $$
and the corresponding regularized determinants, 
$$ {\det_\infty}_{a,\sD}(s) := \exp \left( -\frac{d}{dz} 
\zeta_{a,\sD}(s,z) |_{z=0} \right). $$

\no Proposition \ref{zetaL} then shows that the
archimedean factor of the Hasse-Weil L-function is given by the 
regularized determinant of a zeta function $\zeta_{a,\sD}$ of the 
spectral triple of Corollary \ref{arith3Cor}, namely the one for 
$a=\sigma^L(w)^2$. One advantage of this point of view is that 
one can now see the archimedean factor of the Hasse-Weil L-function
as an element in the family ${\det_\infty}_{a,\sD}$ associated to 
the noncommutative geometry $(\cA,\sH^\cdot, \sD)$.

\medskip

\no Different representations of the Lie algebra $\cg=\gsl(2)$, for
different choice of the K\"ahler class and the corresponding Lefschetz
operators, were considered also in \cite{Looj}. It would be
interesting to see if one can use this formalism of spectral triples in
that context to further investigate the structure of the resulting
K\"ahler Lie algebra (or of the Neron--Severi Lie algebra of projective
varieties considered in \cite{Looj}).

\medskip

\no In the special case of arithmetic surfaces considered in
\cite{CM}, where $X$ is a compact Riemann surface, the result of
Corollary \ref{arith3Cor} can be related to the ``arithmetic
spectral triple'' of \cite{CM}.

\no Consider the complex ${\rm
Cone}(N)^\cdot =\T^\cdot \oplus
\T^\cdot [+1]$ with differential
$$ d_{Cone} = \left( \begin{array}{cc} \delta & N \\ 0 & -\delta
\end{array}\right). $$

\no Proposition 2.23 of \cite{CM} and \S 4 of \cite{KC} show that,
for $X$ a compact Riemann surface, we have
$$ \H^\cdot({\rm Cone}(N)^\cdot,d_{Cone})
\simeq \H^\cdot(\Ker(N)^\cdot,d') \oplus
\H^{\cdot+1}(\Coker(N)^\cdot,d''). $$
This is isomorphic to
$\H^\cdot(\T^\cdot_u,\delta_u)|_{u=0}$. Moreover, under this
identification, the operator $\Phi$ on the cohomology $\H^\cdot({\rm
Cone}(N)^\cdot,d_{Cone})$ considered in \cite{CM} agrees with the operator
$\Phi$ on $\H^\cdot(\T^\cdot_u,\delta_u)|_{u=0}$.

\section{Analogies with loop space geometry}

\smallskip

\no Besides the original motivating analogy with the case of a
geometric degeneration and the resolution \eqref{nearbyresol} of
the complex of nearby cycles, the cohomology theory at arithmetic
infinity defined by the complex $(\T^\cdot, \delta)$ also bears
some interesting formal analogies with Givental's homological
geometry on the loop space of a K\"ahler manifold (\cf
\cite{Giv}).

\smallskip

\no Let $X$ be a compact K\"ahler manifold, with the symplectic form
$\omega$ representing an integral class in cohomology, such that the
morphism $\omega: \pi_2(X) \to \Z$ is onto.
Let $LX$ denote the space of contractible loops on $X$, and
$\widetilde{LX}$ the cyclic cover with
group of deck transformations
\begin{equation}\label{deck}
\pi_2(X) / \Ker \{ \omega: \pi_2(X) \to \Z \} \cong \Z,
\end{equation}
and with the $S^1$-action that rotates loops.
This covering makes the action functional
\begin{equation}\label{Action}
\A(\phi) = \int_{\Delta} \phi^* \omega, \ \ \ \ \  \forall \phi\in
\widetilde{LX}
\end{equation}
single valued, for $\partial\Delta=S^1$. 
In fact, if $\gamma$ denotes the generator of
\eqref{deck}, we have $\gamma^* \A = \A + 1$. The critical
manifold $Crit(\A)=Fix(S^1)$ consists of a trivial cyclic cover of
the submanifold of constant loops $X$. Formally, one can consider
an equivariant Floer complex for the functional $\A$, which is the
complex \eqref{alghbar} with a differential $d_{S^1}\pm \pi_{*}
\pi^*$, which combines the equivariant differential on each
component of $Crit(\A)$ with a pullback--pushforward along
gradient flow lines of $\A$ between different components of the
critical manifold.

\smallskip

\no More precisely, a {\em formal} setting for the construction of
equivariant Floer cohomologies can be described as follows. On a
configuration space $\sC$, which is an infinite dimensional
manifold with an $S^1$ action, consider an $S^1$-invariant
functional $\A: \sC \to \R$, satisfying the following assumptions:
{\em (i)} The critical point equation $\nabla \A =0$ cuts out a finite
dimensional smooth compact $S^1$-manifold ${\rm Crit}(\A)\subset
\sC$. {\em (ii)} 
The Hessian $H(\A)$ on ${\rm Crit}(\A)$ is non-degenerate in the
normal directions. {\em (iii)} For any $x,y\in {\rm Crit}(\A)$, there is a
well defined locally constant relative index $\ind(x)-\ind(y)\in \Z$.
{\em (iv)} For any two $S^1$-orbits $\sO^\pm$ in ${\rm Crit}(\A)$,
the set $\M(\sO^+,\sO^-)$ of solutions to the flow equation
\begin{equation}\label{flow}
\frac{d}{dt} u(t) + \nabla \A(u(t)) =0, \ \ \ \  \lim_{t\to
\pm\infty} u(t) \in \sO^\pm,
\end{equation}
modulo reparameterizations by translations, is either empty or a
smooth manifold of dimension
$\ind(\sO^+)-\ind(\sO^-) + \dim \sO^+ -1$.
{\em (v)} The manifolds $\M(\sO^+,\sO^-)$
admit a compactification to smooth manifolds with corners, with
codimension one boundary strata
$ \cup_{\ind(\sO^+) \geq \ind(\sO) \geq
\ind(\sO^-)}\M(\sO^+,\sO)\times_{\sO} \M(\sO,\sO^-)$,
and with compatible endpoint fibrations $\pi_\pm :
\M(\sO^+,\sO^-)\to \sO^\pm$.

\smallskip
\no For each component $\sO\subset {\rm Crit}(\A)$ one can then
consider the $S^1$-equivariant de Rham complex
\begin{equation}\label{S1dR}
\Omega^{\cdot,\infty}_{S^1}(\sO) := \Omega^\cdot_{inv}(\sO)\otimes
\C [U,U^{-1}],
\end{equation}
where $U$ is of degree two, so that the total degree of $\alpha
\otimes U^r$, with $\alpha\in \Omega^m_{\sO}$ is $i=m+2r$. The
differential is of the form
\begin{equation}\label{dS1}
 d_{S^1}( \alpha \otimes U^r ) = d\alpha \otimes U^r +\iota_V(\alpha)
\otimes U^{r+1}, \ \ \  d_{S^1}U=0,
\end{equation}
where $\iota_V$ denotes contraction with the vector field
generated to the $S^1$-action on $\sO$. The complex \eqref{S1dR}
computes the periodic $S^1$-equivariant cohomology of $\sO$,
$$ H^\cdot(\Omega^{\cdot,\infty}_{S^1}(\sO), d_{S^1}) =
H^\cdot_{S^1,per}(\sO;\C),
$$
which is the localization of $H_{S^1}^\cdot(\sO;\C)$ obtained by
inverting $U$. The Floer complex is then defined as
\begin{equation}\label{S1Floer}
CF^{\ell,\infty}_{S^1} = \oplus_{\ell=\ind(\sO) + m+2r}\,\,
\Omega^{m+2r, \infty}_{S^1}(\sO),
\end{equation}
with the relative index $\ind(\sO)$ computed with respect to a
fixed base point in ${\rm Crit}(\A)$, and with the Floer
differential given by
\begin{equation}\label{FloerD}
D_{\sO^+,\sO^-} (\alpha\otimes U^r) = \left\{ \begin{array}{ll}
d_{S^1}\,
(\alpha\otimes U^r) & \sO^+=\sO^- \\
(-1)^{m} (\pi_{+\,*} \pi_-^* \, \alpha) \otimes U^r &
\ind(\sO^+)\geq \ind(\sO^-) \\
0 & \text{ otherwise, }
\end{array}\right.
\end{equation}
where $\pi_\pm : \M(\sO^+,\sO^-)\to \sO^\pm$ are the endpoint
projections. The (periodic) equivariant Floer cohomology is the
$\C[U,U^{-1}]$ module
\begin{equation}\label{HFinf}
HF_{S^1}^{\cdot ,\, \infty}(\sC;\A) := H^\cdot(CF^{\infty}_{S^1},
D).
\end{equation}
The property $D^2=0$ for the Floer differential holds because of
the structure of the compactification of the
spaces $\M(\sO^+,\sO^-)$ and its compatibility with the endpoint
fibrations (\cf \cite{AB}). The periodic equivariant Floer
cohomology is related to equivariant Floer cohomology and homology
via a natural exact sequence of complexes (\cf \cite{MW2}), of the
form
\begin{equation}\label{splitpm}
 0 \to CF^{*,+}_{S^1} \to CF^{*,\infty}_{S^1} \to CF^{*,-}_{S^1} \to
0,
\end{equation}
where $CF^{*,+}_{S^1}$ is defined as in \eqref{S1Floer}, but with
$\C[U]$ instead of $\C[U,U^{-1}]$ in \eqref{S1dR}. The equivariant
Floer cohomology is defined as
\begin{equation}\label{HF+}
 HF^*_{S^1}(\sC;\A) = HF^{*,+}_{S^1}(\sC;\A)=H^*(CF^{*,+}_{S^1},D),
\end{equation}
while the quotient complex $CF^{*,-}_{S^1}$, with the induced
Floer differential $D^-$, computes the equivariant Floer homology
\begin{equation}\label{HF-}
 HF_{*,S^1}(\sC;-\A) = HF^{*,-}_{S^1}(\sC;\A).
\end{equation}

\smallskip

\no There is also, in Floer theory, an analog of the weight filtration
$W_\cdot$ given by the increasing filtration of the complex
\eqref{S1Floer} by index of critical orbits, $\ind(\sO)\geq k$,
\begin{equation}\label{filtr}
 W_k CF^{\ell,\infty}(\sC,\A) = \oplus_{\ind(\sO)\geq k,\,\,
 i+\ind(\sO)=\ell}\,\, \Omega^i_{S^1}(\sO).
\end{equation}
In the cases where the boundary components corresponding to flow
lines in $\M(\sO^+,\sO^-)$ vanish, the exact sequence collapses
and the Floer cohomology is the equivariant cohomology of the
critical set, \cf \cite{Giv} and \cite{GK}. We have then an
analog, in this context, of the Picard--Lefschetz filtration in
the form \eqref{monof}, by setting
\begin{equation}\label{Afloer}
AF^{s,k}:= CF^{s+k+1,\infty}_{S^1}/ W_k CF^{s+k+1,\infty}_{S^1},
\end{equation}
as the analog of \eqref{nearbyresol}, and
\begin{equation}\label{PLfloer}
L_\ell AF^{s,k} = W_{2k+\ell+1} CF^{s+k+1,\infty}_{S^1}/ W_k
CF^{s+k+1,\infty}_{S^1},
\end{equation}
with
\begin{equation}\label{grLWfloer}
gr^L_\ell AF^{s,k} =\left\{ \begin{array}{ll}
gr^W_{2k+\ell+1}CF^{s+k+1,\infty}_{S^1} & \ell+k \ge 0 \\[2mm]
0 & \ell+k<0 \end{array}\right.
\end{equation}
where
$$ gr^W_{2k+\ell+1}CF^{s+k+1,\infty}_{S^1}=
\oplus_{\ind(\sO)=2k+\ell+1, s-i=\ell+k} \Omega^i_{S^1}(\sO). $$

\smallskip

\no In the case of the loop space of a smooth compact symplectic
manifold $X$, the action functional \eqref{Action} is degenerate,
the transversality conditions fail and the setup of Floer theory
becomes more delicate \cf \cite{PSS}. The argument of \cite{GK}
shows that, in the case of the loop space of a K\"ahler manifold,
the components $\M(\sO^+,\sO^-)$ contribute trivially to the Floer
differential, hence the equivariant Floer cohomology is computed
by the $E_1$ term of the spectral sequence associated to the
filtration $W_\cdot$, namely by the infinite dimensional vector
space $H^\cdot (X;\C)\otimes \C[U,U^{-1}] \otimes
\C[\hbar,\hbar^{-1}]$, where $\hbar$ implements the action of $\Z$
in \eqref{deck}.

\smallskip

\no Thus, one can see a formal analogy between the complex $\fC^\cdot$
of \eqref{cutoffgammaC} and an equivariant Floer complex on the loop
space. The cutoff $k\geq 0$ is then interpreted as a filtration by 
sublevel sets, corresponding to considering $HF^\cdot_{S^1}$ of
$\sC_0=\cA^{-1}(\R_{\geq 0})$. The cutoff $2r+m\geq 0$ instead
corresponds in this analogy to a cutoff on the total degree of the 
equivariant differential forms, $HF^{\cdot\geq 0}_{S^1}$. Finally, the
cutoff $k\geq r+q$ in \eqref{lambda-cut} can be compared to a
splitting of the form \eqref{splitpm}, adapted to the Hodge
filtration. As in the case of the analogy with the resolution of the complex
of nearby cycles, also in this analogy with Floer theory on loop
spaces there is however a fundamental difference in the differentials.
In fact, the term $d_{S^1}$ of the Floer
differential is replaced in the theory at
arithmetic infinity by $\hbar d +d''$, which would not give a
degree one differential on the Floer complex (except in special
cases, like hyperk\"ahler manifolds, where all components of ${\rm
Crit}(\A)$ have relative index zero).

\smallskip

\no By the results of Givental, the equivariant Floer cohomology of
the loop space also has an action of the ring $\sD$ of
differential operators on $\C^*$, where $Q$ acts as
$\hbar=\gamma^*$ and $P$ as a combination of the symplectic form
and action functional on the loop space, \cite{Giv}. This
structure on the Floer cohomology plays an important role in the
phenomenon of mirror symmetry. It is interesting to remark that 
there is a conjectural mirror
relation between the monodromy and the Lefschetz operators (\cf
\cite{GLO} \cite{Manin}). Thus, the question of developing a more
precise relation between the complex $(\T^\cdot,\delta)$ and Floer
theory, with the actions \eqref{actionD} of the ring of
differential operators on $\C^*$, may be interesting in this respect, 
in view of the possibility of addressing such mirror symmetry
questions in the context of arithmetic geometry, by adapting to
arithmetic cohomological constructions the setting of 
homological geometry on loop spaces.

\end{document}